\documentclass[11pt]{article}
 \usepackage{latexsym}
 \usepackage{amsbsy}
 \usepackage{amssymb}
 \usepackage{amsmath}
 \input{amssym.def}
 \input{amssym}
 \usepackage{graphicx}
 \usepackage{subcaption}
 \usepackage{epsfig}
 \usepackage[latin1]{inputenc}%needed for bibliography
 \usepackage{a4wide}
 \usepackage{caption}
 \usepackage{comment}
 \usepackage{epstopdf}
 
 \usepackage{marginnote}
 \usepackage{xcolor}
\usepackage{lineno}

 \newtheorem{thm}{Theorem}[section]
 \newtheorem{defin}[thm]{Definition}
 \newtheorem{lem}[thm]{Lemma}
 \newtheorem{prop}[thm]{Proposition}
 \newtheorem{cor}[thm]{Corollary}
 \newtheorem{rem}[thm]{Remark}
 \newtheorem{ex}[thm]{Example}

 \newcommand{\bthm}{\begin{thm}}
 \newcommand{\ethm}{\end{thm}}
 \newcommand{\bd}{\begin{defin}}
 \newcommand{\ed}{\end{defin}}
 \newcommand{\blem}{\begin{lem}}
 \newcommand{\elem}{\end{lem}}
 \newcommand{\bcor}{\begin{cor}}
 \newcommand{\ecor}{\end{cor}}
 \newcommand{\bprop}{\begin{prop}}
 \newcommand{\eprop}{\end{prop}}
 \newcommand{\brem}{\begin{rem} \rm}
 \newcommand{\erem}{\end{rem}}
 \newcommand{\bex}{\begin{ex} \rm}
 \newcommand{\eex}{\end{ex}}

 \newcommand{\pr}{\noindent{\bf Proof. }}
 \newcommand{\ep}{\nolinebreak{\hspace*{\fill}$\Box$ \vspace*{0.25cm}}}

 \newcommand{\beq}{\begin{equation}}
 \newcommand{\eeq}{\end{equation}}

 \newcommand{\bea}{\begin{eqnarray}}
 \newcommand{\eea}{\end{eqnarray}}

 \newcommand{\beas}{\begin{eqnarray*}}
 \newcommand{\eeas}{\end{eqnarray*}}

 \newcommand{\beqs}{\begin{equation*}}
 \newcommand{\eeqs}{\end{equation*}}

 \newcommand{\bi}{\begin{itemize}}
 \newcommand{\ei}{\end{itemize}}

 \newcommand{\ben}{\begin{enumerate}}
 \newcommand{\een}{\end{enumerate}}

 \newcommand{\ba}{\begin{array}}
 \newcommand{\ea}{\end{array}}

 \newcommand{\R}{\mathbb R}
 \newcommand{\N}{\mathbb N}
 \newcommand{\C}{\mathbb C}

 \newcommand{\cC}{\ensuremath{{\cal C}}}

 \newcommand{\cT}{\ensuremath{{\cal T}}}

 \newcommand{\eps}{\varepsilon}
 \newcommand{\vphi}{\varphi}

 \newcommand{\pd}{\ensuremath{\partial}}

\newcommand{\pa}{\partial}

 \def\dj{d\kern-0.4em\char"16\kern-0.1em}
 \def\Dj{\mbox{\raise0.3ex\hbox{-}\kern-0.4em D}}

\begin{document}
%%%%%%%%%%%%%%%%%
%%%%%%%%%%%%%%%%%

 \title{Global Controllability for Quasilinear Non-negative Definite System of ODEs and SDEs}

 \author{Jasmina Djordjevic\footnote{Faculty of Mathematics and Natural Sciences, University of Nis, Visegradska 33, 18000 Nis, Serbia, Electronic mail: nina19@pmf.ni.ac.rs} \\
 	   Sanja Konjik\footnote{Faculty of Sciences, University of Novi Sad, Trg D.Obradovi\'ca 4, 21000 Novi Sad, Serbia,
         Electronic mail: sanja.konjik@dmi.uns.ac.rs}\\
         Darko Mitrovi\'c
         \footnote{Faculty of Mathematics, University of Vienna, Oskar Morgenstern Platz-1, 1090 Wien, Austria,
         Electronic mail: darko.mitrovic@univie.ac.at}
      \\ Andrej Novak
         \footnote{ Department of Physics, Faculty of Science, Bijeni\v cka cesta 32, University of Zagreb, Croatia,
         Electronic mail: andrej.novak@phy.hr\newline A.Novak is the corresponding author.}}

 \date{}
 \maketitle

 \begin{abstract}
We consider exact and averaged control problem for a system of quasi-linear ODEs and SDEs with a non-negative definite symmetric  matrix of the system.  The strategy of the proof is the standard linearization of the system by fixing the function appearing in the nonlinear part of the system, and then applying the Leray-Schauder fixed point theorem. We shall also need the continuous induction arguments to prolong the control to the final state which is a novel approach in the field. This enables us to obtain controllability for arbitrarily large initial data (so called global controllability). %We also provide a criterion for controllability of the considered systems. 

 \vskip5pt
 \noindent
 {\bf MSC (2010):}
 %\begin{classcode}
  Primary: 34H05; Secondary: 49J15, 93C15, 60H10 \\
 %\end{classcode}
 \end{abstract}

 %Communicated by Firdaus E. Udwadia 
%\nocite{*}
%%%%%%%
\section{Introduction}
%\label{sec:intro}
%\linenumbers
A fundamental issue in control theory is how to choose a source term (usually called the control) in the  equations of the particular model which would govern the system from the given initial state to the prescribed final state. Such results clearly have a wide range of possible applications which makes control theory one of the most attractive fields in mathematics. It is not easy to say where the origins of the substantial mathematical treatment of the theory lie, but one can find lots of information in standard books such as \cite{Lions, book, Zua-11}. Concerning the direction in the control theory that we are going to pursue here, we mention \cite{LZ, LoZ, Pi3, Zua-22} and the references therein.

%As for the recently introduced averaged control which, in particular, we will deal with here, in the case of ODEs the notion is introduced in \cite{Zua-22} and extended in \cite{LZ} on PDEs governed systems.

%In the stochastic situation, subject to the stochastic forcing as in many applications, one cannot have complete and precise grasp over the monitored process.

%The theory is well developed for the systems governed by linear ODEs, but, interestingly enough, unlike the optimal control theory (dealing with minimizing a cost functional; see e.g. \cite{CBG}, references therein and the subsequent paper citations), nonlinear situations have been rarely considered. We mention results in the frame of scalar conservation laws (e.g. \cite{ancona, boria}) where the authors investigate the attainable set of function using the method of generalized characteristics. In other words, their methods of proof are related to the analysis 
%of system of characteristics i.e. system of nonlinear ODEs, but the corresponding initial data play the role of the control and such methods are not applicable here.     

The theory is well developed for both linear and non-linear systems. The theory started with linear systems and most of the related questions have been answered (see e.g. \cite{book}). Nonlinear systems have been intensively considered recently and one can find various problems of this type in \cite{BS, CU, LZh, Naito, UK, U, XZh} and references therein. Let us briefly recall the related results from the field. 

The most general problem in this direction is the nonlinear system of the form
\begin{equation}
\label{general}
\dot{y}(t)=f(t,y(t),u(t)), \ \ y(0)=y^0 \in \R^d,
\end{equation} and the aim is to find the vector valued function $u$ which steers the system to a given final state $y(T)=y^T$. This system was considered in \cite{Naito} and it has been proved that the system is locally approximate null controllable provided that the linear operator of the system is approximately invertible and the linear approximation to \eqref{general} is locally null controllable. A similar problem -- a quasilinear delay integrodifferential equation has been considered e.g. in \cite{BS} (similar problems can be found also in \cite{N1, NP} and many others chronologically subsequent papers). 

The proofs of controllability of nonlinear systems essentially rely on the linear theory and Leray-Schauder-Tikhonov fixed point type theorems. Roughly speaking, we consider the Taylor expansion of $f$ in the neighborhood of $(t,y,u)=(t,0,0)$:
$$
f(t,y,u)=f(t,0,0)+f_y(t,0,0)y+f_u(t,0,0)u+g(t,y,u),
$$ and then we linearize \eqref{general} by considering 
\begin{equation}
\label{gen-lin}
\begin{split}
&\dot{y}(t)=f(t,0,0)+f_y(t,0,0)y+f_u(t,0,0)u+g(t,z,u), \\ &y(0)=y^0 \in \R^d, \ \ y(T)=y^T\in \R^d,
\end{split}
\end{equation} for a fixed function $z:[0,T] \to \R^d$ of appropriate regularity. One proves or assumes existence of the control $u$ such that \eqref{gen-lin} holds and thus one obtains the mapping 
$$
{\cal T}(z)=y
$$  which adjoins the state $y$ to the previously fixed function $z$. Clearly, if we prove existence of a fixed point for the mapping ${\cal T}$ we will prove existence of control to \eqref{general}. 

In the current contribution, we shall consider an ODE system with quasilinear non-negative definite symmetric right-hand side (sometimes called degenerate parabolic system; see \eqref{GenForm}). Such and similar type of equations have been intensively considered recently in both deterministic and stochastic setting (see \cite{BKN} and references therein). In particular, in \cite{XZh} the control problem of the following form has been investigated
\begin{equation}
\label{general-1}
\dot{y}(t)=-A(t, y(t)) y(t)+B(t,y(t))u(t), \ \ y(0)=y^0, \ \ y(T)=y^T, 
\end{equation} where $(t,y) \mapsto A(t,y) \in \R^{d\times d}$ and $(t,y) \mapsto B(t,y) \in \R^{d\times n}$ are matrix-valued functions of appropriate regularity. 
In that case, the situation is more straight forward because we can avoid the Taylor expansion and simply replace $y$ appearing in $A$ and $B$ by $z$ and then apply linear theory to get the state ${\cal T}(z):=y$ and the control $u$ satisfying
\begin{equation}
\label{gen-lin-1}
\dot{y}(t)=-A(t, z(t)) y(t)+B(t,z(t))u(t), \ \ y(0)=y^0, \ \ y(T)=y^T. 
\end{equation} As before, we look for the fixed point of the mapping ${\cal T}$ which provides solution to \eqref{general-1}.

The first obstacle that we encounter is the controllability problem formulated in \eqref{gen-lin-1}. In general, it is not simple to prove that for any $z$ one can find $y$ and $u$ such that \eqref{gen-lin-1} holds, and it is often simply taken as an assumption (see e.g. \cite[(iii)]{BS}). However, in special situations such as the one considered in \cite{LZh}, one can prove the controllability for any $z$ of appropriate regularity. The extenuating circumstance in that case is the assumption that the equation is of strictly parabolic type i.e. the corresponding operator is coercive, and this enables one to obtain the result of such a generality.
The coercivity of the operator then provides an observability result which is, as well known \cite[Ch III]{book}, equivalent to the controllability. 

We also note that the result from \cite{LZh} holds under the assumption that the initial data are small enough (it is called local controllability) and that here, we introduce a method for overcoming the latter confinement. More precisely, by introducing a change of variables $\tau=Kt$ for $K$ large enough, we are able to obtain estimates that provide compactness of the mapping ${\cal T}$ in a small interval $[0,T/K]$ around zero and thus existence of the fixed point in that interval. By continuing the procedure beginning from $T/K$ instead of $0$, and then marching forward step by step, we reach the final moment $T$. 

The paper is organized as follows. In the next section, we introduce necessary notions and notations and formulate precisely the results that we are going to prove. We also provide additional historical facts. In Section 3, we prove that problem \eqref{general-1} admits a solution (see Theorem \ref{main-T}), while in Section 4, we prove the stochastic variant of the result (see Theorem \ref{main-S}). In Section 5, we provide an application in the frame of the population dynamics.  The control function represents a source which means that by choosing appropriate $u$ in \eqref{GenCond} we choose the rate at which we add new individuals into the system or we simply change living conditions so that the population increases or decreases to the desired level. We use an iterative method combined with adaptation of final data in order to cancel computational errors and obtain satisfactory simulations.

\section{Notions, notations, and formulation of the results}

Let us now formally introduce the system that we are going to deal with. To this end, we need the complete probability space
$(\Omega, {\cal F},  {\cal F}_t, {\bf P})$ with the sample space $\Omega$, the $\sigma$-algebra ${\cal F}$, the natural filtration  $\{ {\cal F}_t\}_{t\in [0,T]}$ generated by the standard $n$-dimensional ($n\geq 1$) Wiener process $W_t(\cdot)$, and the probability measure ${\bf P}$.  For an Euclidean space $H$ (such as $\R^d, \R^n,\R^{d\times n}$, $\R^{d\times d}$ depending on the situation), we denote the spaces:
\begin{align*}
L_{{\bf P}}^2(\Omega;C([0,T]);H)&=\{ y: [0,T]\times \Omega\to H | \ \  \forall \omega\in \Omega, \\ & \qquad \qquad y(\omega,\cdot)\in C([0,T];H), \, {\rm and} \, \int_\Omega \sup\limits_{t\in [0,T]}|y(t,\omega)|^2 d{\bf P}(\omega) <\infty\},\\
L_{{\bf P}}^2(\Omega;L^2([0,T]);H)&=\{ y: [0,T]\times \Omega\to H | \, \int_\Omega \int_{[0,T]}|y(t,\omega)|^2 dt d{\bf P}(\omega) <\infty\}.
\end{align*} 
To simplify the notation, in the sequel we imply
$$L^2(\Omega;C([0,T])):=L_{{\bf P}}^2(\Omega;C([0,T]);H), \qquad L^2(\Omega;L^2([0,T])):=L_{{\bf P}}^2(\Omega;L^2([0,T]);H).$$
The most general form of the system to be considered here is
\begin{align}
\label{GenForm} 
dy(t,\omega)=(-A(\omega,t, y(t,\omega)) y(t,\omega)+B(\omega,t,y(t,\omega))u(t,\omega))dt+ Z(t,\omega) dW(t),
\end{align}
where 
\begin{itemize}
\item $y(t,\omega): [0,T]\times\Omega\to \R^{d}$, $y\in L^2(\Omega; C([0,T]))$, is adapted with respect to the filtration ${\cal F}_t$;
\item $Z(t,\omega): [0,T]\times\Omega \to \R^{d\times n}$, $Z\in L^2(\Omega; L^2([0,T]))$, is adapted with respect to the filtration ${\cal F}_t$; 
\item $u(t,\omega): [0,T]\times\Omega\to \R^{n}$, $u\in L^2(\Omega; L^2([0,T]))$, is adapted with respect to the filtration ${\cal F}_t$.
\end{itemize} The function $y$ is usually called the state of the system, while $u$ is the control. The pair $(y,Z)$ is uniquely determined via the corresponding backward differential equation (see \cite{moi1}) which means that we are actually looking only for $y$ and $u$. Let us remark that this is both natural and intuitive  since $Z$ is random and cannot be controlled by initial and/or final data. 

Now we are ready to formally introduce definitions of controllability that we shall need here.

\begin{defin}  We say that the  system \eqref{GenForm} is {\it exactly controllable} if we can find the function $u\in C([0,T];\R^n)$ such that solution $y\in C^1([0,T];\R^d)$ to \eqref{GenForm} satisfies  
\begin{equation*}
y(0)=y^0,    \qquad y(T)=y^T
\end{equation*} for any $y^0,y^T\in \R^d$.
 \end{defin}\label{exa.con} We note that we shall prove the exact controllability in the deterministic case i.e. in the case when \eqref{GenForm} does not depend on $\omega\in \Omega$.

\begin{defin} We say that the system \eqref{GenForm} admits the averaged control if we can find $u \in L^2(\Omega;L^2[0,T])$, such that $(y,Z)$, $y\in L^2(\Omega;C([0,T]))$, $Z\in L^2(\Omega;L^2[0,T])$, solves \eqref {GenForm}, and the following initial and final conditions are satisfied:
\begin{equation}
\label{GenCond}
y(0)=y^0,    \qquad E(y(T,\cdot))=y^T
\end{equation} where $y^0$ is an $d$-dimensional random variable, $y^T$ is a constant vector  ($E$ is the expectation operator).  In this case system (\ref{GenForm}),  (\ref{GenCond}), is called {\it controllable in average}.  \label{ave.con}
\end{defin}

Although we consider a system of ODEs, the techniques used here can be transferred to parabolic PDEs via the Galerkin approximation (this will be a subject of future research). Lack of strict parabolicity is a major obstacle in the proof of the solvability of \eqref{gen-lin} so we shall need to assume it below. However, in some important special situations (see \cite{LZh} and Proposition \ref{p1.1}), we are able to overcome this obstacle. More importantly, unlike the situation from \cite{LZh}, we are able to prove the controllability result for \eqref{GenForm} for arbitrary large initial data of appropriate regularity.  

%To this end, it is necessary to get appropriate compactness for the mapping $T$
Let us now proceed with the assumptions that we shall need. We shall imply:
\begin{itemize} 

\item[(i)] $A:\Omega \times \R^+ \times \R^d\to \R^{d\times d}$ is a smooth non-negative definite symmetric matrix-valued process; 

%\item[(i')] $A(y,\omega):\R^d\times \Omega\to \R^{d\times d}$ is a smooth bounded matrix-valued mapping;

\item[(ii)] $B: \Omega\times \R^+ \times \R^d\  \to \R^{d\times n}$ is a matrix-valued $L^\infty$-uniformly bounded process; 

\item[(iii)] For every $v=(v_1,\dots,v_d)$, $v_j\in L^2(\Omega;C([0,T]))$, $j=1,\dots,d$, the Gramian (note the difference in the notation of the final time $T$ and the transpose ${\rm T}$ below)
\begin{equation}
\label{gram}
G_c(0,T)=E\left(\int_0^T e^{\int_0^t -A(\omega,t,v) dt'} B(\omega,t,v) B(\omega,t,v)^{\rm T} \left( e^{\int_0^t(-A(\omega,t,v))dt'}\right)^{\rm T} dt\right)
\end{equation} is invertible.
\end{itemize} In the deterministic case we have $Z\equiv 0$ and we do not have dependence on the variable $\omega$. Condition (iii) thus becomes:

\begin{itemize}
\item[(iii)'] If in \eqref{GenForm} we have $Z\equiv 0$ and the coefficients are independent of $\omega\in \Omega$, then we assume
\begin{equation}
\label{gram-determ}
G_c(0,T)=\int_0^T e^{\int_0^t -A(t,v) dt'} B(t,v) B(t,v)^{\rm T} \left( e^{\int_0^t(-A(t,v))dt'}\right)^{\rm T} dt
\end{equation} is invertible for every $v\in C([0,T])$. 

\end{itemize} We note that the invertibility of the Gramian corresponding to a linear controllability system is necessary and sufficient for the controllability to hold (see \cite[Theorem 3.2]{4C} in the stochastic and \cite{book} in the deterministic case). Here, urged by the nonlinearity of the system and motivated by the construction of the proof, we needed to modify the condition in the form given in (iii).  
Although condition (iii) seems too strong in the sense that it is highly non-trivial to verify, it is a standard fare in this kind of research (see again \cite[(iii)]{BS}). However, in some interesting and important special cases, one can simplify the condition. In \cite{LZh}, the controllability has been proved for strictly parabolic equations (even partial differential equations). Here, we shall prove controllability for a porous media type equation (which is an example of a degenerate parabolic equation; see \cite{stoch7, BPR}):
\begin{equation}
\label{porous-med} 
y'(t)=-A|y(t)|^m y(t)+Bu(t), \qquad y(0)=y^0,
\end{equation} where $A\in \R^{d\times d}$ and $B\in \R^{d\times n}$ are constant matrices. We shall pay special attention to the equation in Section 5 where we shall introduce a numerical procedure and provide corresponding computer simulations. 

Interestingly enough, controllability for \eqref{porous-med} is actually equivalent to the Kalman rank condition. Let us prove this fact before we continue with the introduction.

\begin{prop}
\label{p1.1}
Let $A\in \R^{d\times d}$ and $B\in \R^{d\times m}$ be constant matrices. If the Kalman rank condition
\begin{equation}
\label{kalman}
{\rm rank}\big[B\big| AB \big| \dots \big| A^{d-1}B \big]=d
\end{equation} holds, then the Gramian corresponding to \eqref{porous-med}
\begin{equation}
\label{gram-d}
G_c(0,T)=\int_0^T e^{\int_0^t(-|v|^m A)dt'} B B^{\rm T} \left( e^{\int_0^t(-|v|^m A )dt'}\right)^{\rm T} dt
\end{equation}
is invertible for every $v\in C([0,T])$ (in other words, (iii') holds).
\end{prop} 
\begin{rem}
Let us note that here we do not require non-negative definitness or symmetricity of the matrix $A$. However, it is necessary for the proof of global controllability of \eqref{porous-med}. 
\end{rem}
\pr
Invertibility of \eqref{gram} is equivalent to the complete controllability of the system
\begin{align}
\label{linearized-I}
y'(t)=-A|v(t)|^m y(t)+Bu(t)\\
\label{init-fin-I} 
y(0)=y^0, \ \ y(T)=y^T
\end{align} for any $y^0, y^T\in \R^d$ (see \cite[p. 77]{book}). In other words, if \eqref{gram-determ} is invertible, then for every $y^0,y^T$ we can find $y$ and $u$ so that \eqref{linearized-I}, \eqref{init-fin-I} are satisfied and vice versa. Thus, to prove the proposition, it is enough to show that there exist $y$ and $u$ such that \eqref{linearized-I}, \eqref{init-fin-I} hold. To this end, we consider the following perturbation of \eqref{linearized-I}
\begin{equation}
\label{pert}
y_\eps'(t)=(-A|v(t)|^m-\eps A) y_\eps(t)+Bu_\eps(t), \ \ \eps>0.
\end{equation}  We shall show that if the Kalman rank condition is satisfied then \eqref{pert} is completely controllable i.e. for every $y^0,y^T$ we can find $y_\eps$ and $u_\eps$ so that \eqref{pert}, \eqref{init-fin-I} are satisfied. Consider the Gramian corresponding to \eqref{pert}:
$$
G_\eps(0,T)=\int_0^T e^{-\int_0^t(|v|^m+\eps) A dt'} B B^{\rm T} \left( e^{-\int_0^t(|v|^m+\eps) A dt'}\right)^{\rm T} dt.
$$
Now, let us prove that invertibility of $G_\eps$ is equivalent to \eqref{kalman}.

By contradiction, assume that $G_\eps(0,T)$ is not invertible. This means that there exists $w\in \R^d\backslash \{0\}$ such that
\begin{equation}
\label{eps1}
\begin{split}
&w\cdot G_\eps(0,T)\cdot w^{\rm T}=0 \implies \int_0^T w\cdot e^{-\int_0^t(|v|^m+\eps) A dt'} B B^{\rm T} \left( e^{-\int_0^t(|v|^m+\eps) A dt'}\right)^{\rm T} \cdot w^{\rm T} dt=0\\
& \implies \int_0^T |w\cdot e^{-\int_0^t(|v|^m+\eps) A dt'} B|^2 dt=0 \implies w\cdot e^{-\int_0^t(|v|^m+\eps) A dt'} B \equiv 0 
\end{split}
\end{equation} for every $t\in [0,T]$. From \eqref{eps1} we get
$$
w\cdot e^{-\int_0^t(|v|^m+\eps) A dt'} B=0 \implies w\cdot \left(I+\sum\limits_{j\in \N} A^j \frac{(-\int_0^t(|v|^m+\eps)dt')^j}{j!}\right) B=0.
$$ Since $\int_0^t(|v|^m+\eps)dt'>0$ and since $t\geq 0$ is arbitrary, we conclude 
$$
w \cdot  A^j B=0, \ \ \forall j\in \N,
$$ implying that the Kalman condition is not satisfied. Let us remark in passing that according to the Cayley--Hamilton theorem, $k>n$-th power of the matrix $A$ can be represented as a linear combination of powers of order $0,1,\dots,n$ and thus the opposite statement holds as well. 

From the above, we see that the Kalman condition implies invertibility of the Gramian $G_\eps(0,T)$ independently of $v$. Thus, we see that \eqref{pert} is completely controllable. Moreover, we know that the state function $y_\eps$ and the control $u_\eps$ are uniformly bounded with respect to $\eps$ (see \cite[p.77]{book}). Thus, the families $(y_\eps)$ and $(u_\eps)$ admit weak-$\star$ limits along a common subsequence which we denote by $y$ and $u$, respectively. Clearly, the (vector) functions $y$ and $u$ satisfy \eqref{linearized-I}, \eqref{init-fin-I} from where invertibility of Gramian \eqref{gram-d} follows. \ep

In the stochastic variant of \eqref{porous-med} (to be considered numerically in Section 5), we aim for the average control \cite{Zua-22} and in this case the Kalman rank condition should be modified as in \cite[Theorem 1]{Zua-22}. The proof, however, goes along the same line as the one of Proposition \ref{p1.1}.  
Moreover, we propose a numerical procedure which seems quite robust and converges very fast (see Section 5) and whose stability can be used as a practical test of the controllability.

Continuing in this direction, the second aim of the paper is to examine controllability of the system \eqref{GenForm} in the sense that we are looking for the control (vector) function $(u,Z)$ which steers the state function $y$ from the given initial state to the given final state in the sense of expectation (see \eqref{GenCond} below). A similar problem has been considered in \cite{gashi, 4C, wang}, but with linear coefficients. We also note that it is possible to include more general SDEs as in \cite{4C}, but the problem is essentially reduced to the linear variant of the system of the form \eqref{GenForm} (see \cite[Theorem 3.4]{4C}).

We conclude this section by noting that the general stochastic case of the described control problem \eqref{GenForm} refers to the theory of backward stochastic differential equations (BSDEs). BSDEs is a stochastic differential equation of the following form
$$dy(t)=f(t,y(t),z(t))dt+z(t)dW(t), \quad  y(T)= \xi, $$
where the function $f$ is the generator and $\xi$ is the final condition. The solution of this equation consists of an ordered pair of adapted processes $(y, z)$ that satisfies this equation. In the field of stochastic control, the adapted process $z(t)$ is considered as a process of control, and the process $y(t)$ is the process of the state of the system. The goal is to define an adapted process $z(t)$ so that the state of the system $y (t)$ is brought to a fixed value of $\xi$ at the time $t = T$. This problem is called the reachability problem. 

Linear type of BSDEs was introduced in \cite{bismut}, and later \cite{Mao,pp1} some results on the existence and uniqueness of the adapted solutions in the nonlinear case are established. Results regarding problems of perturbations for different type of BSDEs, under several conditions are given in \cite{moi1,moi2,moi3}. Interestingly, the BSDE given above in general uniquely determines the functions $y$ and $z$. Therefore, if we add the initial condition as in \eqref{GenForm}, we need another function (the function $u$) to control the process.

 \section{A control problem for nonlinear system of ODEs}
 %\label{sec:main}
%%%%%%%

 We are concerned with a control problem for the nonlinear system of ODEs of the following form
 \bea
 \frac{dy}{dt}(t) &=& - A(t,y(t)) y(t) + B(t,y(t)) u(t), \quad t\in[0,T], \label{Dproblem} \\ 
 y(0) &=& y^0, \label{Dpuslov} 
 \eea  
 where $A: [0,T]\times \R \to \R^{d\times d}$ is a smooth positive semi-definite matrix-valued mapping, $B: [0,T]\times \R \to \R^{d\times n}$ 
 is a smooth matrix-valued mapping uniformly bounded with respect to $L^\infty$-norm, and $y^0\in\R^d$ is given initial state. 
 The problem consists of determining a function $u\in L^2([0,T]; \R^n)$ such that for the given final state
 $y^T\in\R^d$ and the final time $T$ 
 \beq \label{finstat}
 y(T) = y^T.
 \eeq 
 
% We remark that for $m=0$ the control problem reduces to the well known linear situation. Let us now prove the main theorem of the section.

 \bthm
 \label{main-T}
 Assume that conditions (i), (ii), and (iii) are satisfied. Then, there exist unique functions $y\in C^1([0,T];\R^d)$ and $u\in C([0,T];\R^n)$ such that the exact control problem \eqref{Dproblem}, \eqref{Dpuslov}, \eqref{finstat} is satisfied.
 \ethm
 
 \pr 
 Our approach in solving this problem will consist of several steps. 
 First, we shall linearize equation \eqref{Dproblem} as follows
 \beq \label{Lproblem}
 \frac{dy}{dt}(t) = - A(t,v(t)) y(t) + B(t,v(t)) u(t), \quad t\in[0,T],
 \eeq
 where $v\in C([0,T])$ is fixed. For such a fixed $v$, we look for the control $u$ steering 
 the system from the initial state $y^0$ to the prescribed final state $y^T$. 

 Since the problem is linear, we can associate to \eqref{Lproblem} an adjoint system in 
 order to use the observability concept and results that are well known for the linear case \cite{book} and even in the stochastic setting \cite[Theorem 3.2]{4C}. Once the linearized problem is solved, we shall employ fixed point arguments in order to reach a solution to the initial problem. 
 
The adjoint system associated to the linearized control problem is given by
 \bea
 -\frac{d\vphi}{dt}(t) &=& - A^*(t,v(t)) \vphi(t), \quad t\in[0,T], \label{vAS1} \\ 
 \vphi(T) &=& \vphi^T, \label{vAS2} 
 \eea where $\vphi^T\in\R^d$ minimizes the convex functional (see \cite{Zuabook})
 \begin{equation}
 \label{CF}
 J(\vphi^T)=\int_0^T |B^*(t,v(t)) \vphi(t)|^2 dt+ \langle y^0 , \varphi(0) \rangle -\langle y^T , \varphi^T\rangle.
\end{equation}  
 
 Then, the so called observability condition reads
 \beq \label{obsL}
 \|\vphi^0\|^2 \leq C \int_0^T \|B^*(t,v(t)) \vphi(t)\|^2 \,dt,
 \eeq  for some constant $C$. It provides solvability of 
 problem \eqref{Lproblem}, \eqref{Dpuslov}, \eqref{finstat}. The standard proof (e.g. \cite{4C} Theorem 4.1.) provides that the 
 controllability of system \eqref{Lproblem} in $[0,T]$ is equivalent to
 observability \eqref{obsL} of system \eqref{vAS1},\eqref{vAS2} on $[0,T]$.  
 Moreover, a solution to problem \eqref{Lproblem} is then given by $u(t)=B^* \vphi(t)$ and we have
 \beq \label{solL}
 \begin{split}
  |u(t)|&=\big| B^*(t,v(t)) \vphi(t)\big|= \big| B^*(t,v(t)) e^{\int_t^T -A(s,v(s))\, ds} \vphi^T \big| \\ 
  & \leq |B^*(t,v(t)) \vphi^T| \leq C  |\vphi^T|,
\end{split} 
 \eeq   where $\vphi$ is a solution to the corresponding adjoint problem. 
 
 It is also well known \cite{book} that condition (iii) from the introduction provides existence of the exact control for \eqref{Lproblem}, \eqref{Dpuslov}, \eqref{finstat}.
 
 Introduce now a change of variables $\tau=K_1 t$ into \eqref{Lproblem}, 
 for some parameter $K_1>0$.
 Denote by 
 \begin{align*}
 &\tilde{v}(\tau)=v\Big(\frac{\tau}{K_1}\Big), \ \ \tilde{u}(\tau)=u\Big(\frac{\tau}{K_1}\Big), \ \ 
 \tilde{y}(\tau)=y\Big(\frac{\tau}{K_1}\Big), \\ 
 &\qquad\qquad \tilde{A}(\tau,\tilde{v}(\tau))=A(\frac{\tau}{K_1},\tilde{v}(\tau)), \ \ \tilde{B}(\tau,\tilde{v}(\tau))=B(\frac{\tau}{K_1},\tilde{v}(\tau)).
 \end{align*}
 Then, using the fact that $\pd_t y=\pd_\tau \tilde{y}\, \frac{\pd\tau}{\pd t}=\pd_\tau \tilde{y} \cdot K_1$,
 we see that \eqref{Lproblem} becomes
 $$
 K_1 \frac{d\tilde{y}}{d\tau}(\tau) = -\tilde{A}(\tau,\tilde{v}(\tau)) \tilde{y}(\tau) + \tilde{B}(\tau,\tilde{v}(\tau))\tilde{u}(\tau), \quad \tau\in[0,K_1 T],
 $$ 
 or, after dividing this equation by $K_1 $, the problem reduces to
 \bea 
 \frac{d\tilde{y}}{d\tau} &=& - \frac{\tilde{A}(\tau,\tilde{v}(\tau))}{K_1}\tilde{y}(\tau) 
 + \tilde{B}(\tau,\tilde{v}(\tau)) \frac{\tilde{u}(\tau)}{K_1}, \quad \tau\in[0,K_1 T], \label{tilS1} \\
 \tilde{y}(0) &=& y^0, \ \ \tilde{y}(K_1 T)=y^T. \label{tilS2}
 \eea 
% Further, $\tilde{y}(K_1  T) = y^F$. In addition we extend the coefficients in \eqref{tilS1} 
% by constant for values $\tau>K_1  T$, so that the continuity of $\tilde{v}$ and $\tilde{u}$ is preserved. 
% This means that the solution $\tilde{y}$ will be globally defined.
 
 Next, introduce an operator $\cT:\cC([0,{T}])\to \cC([0,{T}])$, which maps $\tilde{v}$ to a solution 
 $\tilde{y}$ of \eqref{tilS1}, \eqref{tilS2} (more precisely to the first component of the solution $(\tilde{y}, \tilde{u})$) i.e.,
 $$
 \cT(\tilde{v})=\tilde{y}.
 $$
 
 Furthermore, the solution of \eqref{tilS1}, \eqref{tilS2} is given by
 \beq \label{solll}
 \tilde{y}(\tau) = \exp\Big(-\int_0^\tau \frac{\tilde{A}(s,\tilde{v}(s))}{K_1} \,ds\Big) \left[
 \int_0^\tau \exp \Big(\int_0^s \frac{\tilde{A}(\theta,\tilde{v}(\theta))}{K_1} \,d\theta\Big) \tilde{B}(s,\tilde{v}(s)) \frac{\tilde{u}(s)}{K_1}\,ds
 + y^0\right].
 \eeq   Since an arbitrary $\tilde{v}\in\cC([0,K_1 T])$ is bounded, it follows from \eqref{solL}
 that $\tilde{u}$ is also bounded. 
 Thus for any $\tilde{v}$ one can choose a constant $K_1\geq 1$ such that 
 \beqs %\label{estL}
 \Big| \frac{\tilde{A}(\cdot,\tilde{v}(\cdot))}{K_1} \Big|\leq 1 \quad \mbox{ and } \quad \Big|\tilde{B}(\cdot,\tilde{v}(\cdot))\frac{\tilde{u}}{K_1}\Big|\leq 1.
 \eeqs  Next, since $A$ is symmetric, it holds 
 $$
 |\exp(A)|=\exp({|A|}) = \exp(\max\limits_{i} \lambda_i)
 $$ where $\lambda_i$ are non-negative eigenvalues of the matrix $A$ (the equality holds for the $2$-norm but since all norms are equivalent in finite dimensional spaces, we will keep writing $|\cdot|$). From here and \eqref{solll}, one can estimate $\tilde{y}$ as follows
 \beq \label{y-0-eps}
 |\tilde{y}(\tau)| \leq \big| \exp\Big(-\int_0^\tau \frac{\tilde{A}(s,\tilde{v}(s))}{K_1}\,ds\Big) \big|
 \left[\int_0^\tau \big|\exp\Big(\int_0^s \frac{\tilde{A}(\theta,\tilde{v}(\theta))}{K_1} \,d\theta\Big) {\bf 1}\big| \,ds
 +|y^0|\right] \leq \tau +|y^0|,
 \eeq 
 where ${\bf 1}=(1,\dots,1) \in \R^d$. This estimate together with \eqref{tilS1} implies
 \begin{equation} \label{estimate}
 |\partial_{\tau} \tilde{y}| \leq  \tau+|y^0| + 1.
 \end{equation} 
 From here, we see that for any $\tilde{T}>0$ (and in particular for $\tilde{T}=T$), the derivative $\partial_{\tau} \tilde{y}$ is uniformly bounded on 
 $[0,\tilde{T}]$ independently of $K_1$. Therefore, by the Arzela-Ascolli theorem,
 the operator $\cT:C([0,T])\to C([0,T])$ is compact. 
 Moreover, the set 
 $$\{x\in \C([0,T]):\, x=\lambda Tx, \lambda\in[0,1]\}
 $$ is bounded due to 
 \eqref{y-0-eps}. The latter implies that the conditions of the 
  Leray-Schauder fixed point theorem \cite{GT} are satisfied for the operator $\cT$ which in turn implies existence of 
 a fixed point $\tilde{y}$ of $\cT$ that satisfies 
 $$
 \frac{d\tilde{y}}{d\tau} = -\frac{\tilde{A}(\tau,\tilde{y}(\tau))}{K_1}\tilde{y}(\tau) + \tilde{B}(\tau,\tilde{y}(\tau))\frac{\tilde{u}(\tau)}{K_1}, 
 \quad \tau \in [0,T].
 $$
 By reintroducing the change of variables $t=\tau/K_1$ we obtain the solution $y_1(t)=\tilde{y}(t/K_1)$ to 
 \eqref{Dproblem} on the interval $[0,T/K_1]$.
 
We now repeat the whole procedure for the problem \eqref{Dproblem} with the initial data
given at $T/K_1 $, i.e., $y(T/K_1 ) = y_1(T/K_1 )$, and the same final state $y(T)=y^T$. We thus obtain the function $y_2$ representing the solution to \eqref{Dproblem} on the interval $[T/K_1 ,T/K_1 +T/K_2 ]$. Continuing the procedure, we obtain a sequence of functions $y_n(t)$ satisfying \eqref{Dproblem} on the intervals $[\sum\limits_{i=0}^{n-1} T/K_i ,\sum\limits_{i=0}^{n-1} T/K_i +T/K_n ]$, where we imply $K_0=\infty$. Next, note that the function 
\begin{equation}
\label{contind}
Y^n(t)=y_j(t), \ \ t\in \big[\sum\limits_{i=0}^{j-1} T/K_i ,\sum\limits_{i=0}^{j} T/K_i \big], \ \ j\in \{1,\dots,n \}
\end{equation} solves equation \eqref{Dproblem} on the interval $[0,\sum\limits_{i=0}^{n} T/K_i ] $. In that way, we obtain a sequence of solutions $Y^j$ to \eqref{Dproblem} defined on intervals $[0,\sum\limits_{i=0}^{j} T/K_i ]$.

It remains to be shown that a solution obtained by this construction will eventually reach the final state 
$y^T$ at $t=T$. We shall use a continuous induction argument \cite{Horm}. To this end, assume that the maximal interval on which we can find a function $Y$ satisfying \eqref{Dproblem} and \eqref{contind} is $[0,T_y]$ with $T_y$ strictly less than $T$.
Then, by repeating once again our procedure with the inital data given at $T_y$, we obtain a solution defined 
at $[T_y, T_y+\frac{T}{K_y }]$ - a contradiction with the maximality assumption. Therefore, we obtain 
a desired solution $y$ to \eqref{Dproblem}, \eqref{Dpuslov}, \eqref{finstat}.

The solution is unique since the function $u$ is defined via the solution of the backward problem \eqref{vAS1}, \eqref{vAS2} (see \cite{pp1}) which uniquely depends on the final data. On the other hand, the final datum $\varphi^T$ is a minimizer of the convex functional \eqref{CF} and therefore it is  unique. 
 \ep
 
 \section{Stochastic nonlinear averaged control problem}

In this section, we shall consider the average control problem (Definition 1.1) for a nonlinear SDE. More precisely, we aim  to find $(y,u,Z)$ adapted with respect to the filtration ${\cal F}_t$ such that following equation and conditions are satisfied,  
 \begin{align}
 \label{S-Dprob}
  {dy} &= \left( - A(\omega, t,  y(t,\omega)) y(t,\omega)+ B(\omega,t,y(t,\omega))u(t,\omega) \right) dt + Z(t,\omega) dW(t), \quad t\in[0,T]
 \\
 \label{S-in-fin}
 \qquad y(0,\omega) &= y^0(\omega), \ \ E( y(T,\cdot))=y^T.
 \end{align}  
%As for the matrices $A(\omega) \in \R^{d \times d}$ and $B(\omega)\in \R^{d\times N}$, we assume that they are ${\bf P}$-measurable. Moreover, we assume that that matrix $A(\omega)$ is semidefinite almost surely. 
 
As in the previous section, we first need existence of control for the linear variant of \eqref{S-Dprob} (under conditions (i), (ii), and (iii)):
\begin{align}
 \label{S-Lprob}
  {dy} &= \left( - A(\omega, t,  v(t,\omega)) y(t,\omega)+ B(\omega,t,v(t,\omega))u(t,\omega) \right) dt + Z(t,\omega) dW(t), \quad t\in[0,T]
 \\
 \label{L-in-fin}
 y(0,\omega) &= y^0(\omega), \ \ E( y(T,\cdot))=y^T,
\end{align} for the stochastic process $v\in L^2(\Omega;C([0,T])$. The result is proved in \cite[Section 3.2]{4C} under the stronger requirement of exact controllability. However, the result which can be found there does not contain an explicit bound for the control $u$ which we need here in order to repeat the arguments from the previous section. Therefore, we briefly recall the proof of the control problem \eqref{S-Lprob}, \eqref{L-in-fin}. 

First, let us show that without loss of generality we can assume that $y^0=0$ in \eqref{S-in-fin}. Indeed, we solve \eqref{S-Lprob} with $u\equiv 0$ to obtain the solution $y_1$ to the (forward) stochastic differential equation \eqref{S-Lprob} augmented with the initial data $y_1(0,\omega)=y^0$. Then, we find the functions $y_2$ (the state function) and $(u,Z)$ (the control) which solve the control problem for \eqref{S-Lprob} with the zero initial data and the final condition $E(y_2(T,\cdot))=y^T-E(y_1(T,\cdot))$. The function $y=y_1+y_2$ and the control $(u,Z)$ then represent solution to \eqref{S-Lprob}, \eqref{L-in-fin}. Thus, as we can see, the control problem actually reduces to dealing with the situation when $y^0=0$. We have the following theorem

 \bthm
\label{T-control-S}
Under conditions (i), (ii) and (iii), and the assumption $y^0=0$, the system \eqref{S-Lprob} is controllable in average (Definition 1.1) in the sense that \eqref{L-in-fin} holds for any $y^T \in \R^m$. Moreover, $y\in L^2(\Omega;C([0,T)])$ and $(u,Z)\in L^2(\Omega;L^2([0,T)])$ satisfying \eqref{S-Lprob}, \eqref{L-in-fin} are unique. 

%and the function $u$ is the solution $\varphi$ to the following Cauchy problem for the backward stochastic differential equation 
\ethm
\pr Existence of the control such that \eqref{S-Lprob}, \eqref{L-in-fin} are satisfied can be found in \cite[Section 3.6]{4C} even in the case of exact controllability. Thus, we need to prove the second part of the theorem. To this end, recall that the controllability in the sense \eqref{S-Lprob}, \eqref{L-in-fin} equivalent to the following observability inequality:
\begin{equation}
\label{observ}
C \int_0^T E[(B^* \varphi)^2+|\tilde{Z}|^2] dt \geq \delta \big{|} \varphi^T \big{|},
\end{equation} where $C,\delta$ are constants, and $(\varphi,\tilde{Z})$ is the unique adaptive solution to the BSDE 
\begin{equation}
\label{BSE}
d\varphi=A^*(\omega,t,v) \varphi dt+ \tilde{Z} dW_t, \ \ \varphi(\omega,T)=\varphi^T,
\end{equation} where $\varphi^T$ minimizes the functional

\beq
\label{min}
J(\varphi^T)=\frac{1}{2}\int_0^T \big| E\left( |B^* \varphi|^2+|\tilde{Z}|^2 \right) \big|^2 dt +\langle \varphi^T, y^T \rangle.
\eeq We remark that the functional is convex and thus the minimal value is unique.

 Moreover, the control $u$ is given by
 $$
 u=B^* \varphi.
 $$

Let us briefly recall the arguments leading to the latter conclusions.

We know that a bounded linear operator ${\bf K}$ between two Banach spaces ${\bf X}$ and ${\bf Y}$ is surjective if and only if the corresponding adjoint operator ${\bf K}^*$ satisfies the coercivity inequality (see \cite[Proposition 4.1]{4C}):
\begin{equation}
\label{coerc}
\begin{split}
&\big{|}{\bf K}^* y^*\big{|} \geq \delta |y^*|_{{\bf Y}^*}, \ \ y^*\in {\bf Y}^*.
\end{split}
\end{equation}

In our situation, we set 
\begin{equation*}
\begin{split}
& {\bf K}: L^2(\Omega\times [0,T]; \R^{d\times n})\to\R^d,\\
&{\bf K}(u,Z)=y^T
\end{split}
\end{equation*} and then aim to prove that ${\bf K}$ is surjective using \eqref{coerc}. To this end, from \eqref{S-Lprob} we infer that 
$$
{\bf K}^*(\varphi^T)=(B^*\varphi,\tilde{Z}(t)).
$$ Indeed, from \eqref{S-Lprob} and \eqref{BSE}, using the It\^o formula, we get:
\begin{align*}
d \langle y, \varphi \rangle-\langle y, d\varphi \rangle- \langle Z, \tilde{Z} \rangle dt=\langle -A y, \varphi \rangle dt+ \langle Bu,\varphi \rangle dt
+\langle Z,\varphi \rangle dW(t).
\end{align*} Taking the expectation of the latter expression and then integrating it over $[0,T]$, we have
\begin{align}
\label{***}
&\langle y^T, \varphi^T \rangle- \langle y^0, \varphi^0 \rangle -\int_0^T \langle y, d\varphi \rangle dt- E\big( \int_0^T \langle (u,Z), (B^*\varphi,\tilde{Z})^T \rangle dt \big) \\
&=\int_0^T E(\langle y, -A^* \varphi \rangle) dt,
\nonumber
\end{align} from where, according to \eqref{BSE} and keeping in mind that $y^0=0$, we conclude
\begin{equation}
\label{EL}
\langle {\bf K}(u,Z), \varphi^T \rangle =\langle y^T, \varphi^T \rangle= \int_0^T E (\langle (u,Z),(B^*\varphi,\tilde{Z})^T \rangle) dt
\end{equation}i.e. ${\bf K}^*(\varphi^T)= (B^*\varphi,\tilde{Z})$. From here, we see that coercivity inequality \eqref{coerc} for our operator ${\bf K}$ is equivalent to \eqref{observ}.

Finally, recall that if the Banach spaces ${\bf X}$ and ${\bf Y}$ are reflexive and the map $x^*\mapsto |x^*|^2_{{\bf X}^*}$ is Fr\'echet differentiable and convex, then  the linear mapping ${\bf K}: {\bf X}\to {\bf Y}$ is surjective if and only the functional 
\begin{equation}
\label{J}
{\cal J}(y^*;y)=\frac{1}{2}\big{|} {\bf K}^*y^*\big{|}^2_{{\bf X}^*}+\langle y, y^*\rangle, \ \ y^*\in {\bf Y}^*,
\end{equation}admits a unique optimal minimum over $y^*$ (see \cite[Proposition 4.1]{4C}). In other words, conditions \eqref{coerc} and \eqref{J} are equivalent.

In our case, we take ${\bf X}= L^2(\Omega\times [0,T])$ and ${\bf Y}=\R^n$. Then, we notice that \eqref{EL} are actually Euler-Lagrange equations for \eqref{J} with $y=y^T$.  Thus, computing the first variation of ${\cal J}$ we see that the control 
$$
u(t,\omega)=B^*\varphi(t,\omega), \ \ Z(t,\omega)=\tilde{Z},
$$ ensures that the final condition $y^T$ is satisfied and it minimizes the functional ${\cal J}$.

This completes the proof.
\ep

We can now adapt the procedure from the proof of Theorem \ref{main-T} to obtain the following result.

\bthm
\label{main-S}
Under the conditions (i),(ii), and (iii), the system \eqref{S-Dprob} is controllable in average in the sense that there exist unique $y\in L^2(\Omega;C([0,T)])$ and $u, Z\in L^2(\Omega;L^2([0,T)])$ such that \eqref{S-in-fin} holds for any $y^0\in L^1(\Omega)$ and $y^T \in \R^d$.
\ethm
\pr 
As in the proof of Theorem \ref{main-T}, we introduce the change of variables $\tau=K_1  t$ into \eqref{S-Dprob} 
 for a constant $K_1>1$ and denote

 \begin{equation*}
% \label{COV}
 \begin{split}
 &\tilde{v}(\tau,\omega)=v\Big(\frac{\tau}{K_1 },\omega\Big), \ \ \tilde{u}(\tau,\omega)=u\Big(\frac{\tau}{K_1},\omega \Big), \\ 
 &\tilde{y}(\tau,\omega)=y\Big(\frac{\tau}{K_1 },\omega\Big), \ \ \tilde{Z}(\tau,\omega)=Z\Big(\frac{\tau}{K_1 },\omega\Big), \ \ \tilde{W}_\tau=W_{K_1 t},\\
 & \tilde{A}(\omega, \tau,  \tilde{v}(\tau,\omega))=A(\omega, K_1 t,  y(K_1 t,\omega)), \ \ \tilde{B}(\omega, \tau,  \tilde{v}(\tau,\omega))=B(\omega, K_1 t,  y(K_1 t,\omega)) 
 \end{split}
  \end{equation*} 
 Then, from \eqref{S-Lprob} we have
\begin{equation}
\label{stoch-K}
 d \tilde{y}(\tau,\omega)=(- \frac{\tilde{A}(\omega, \tau,  \tilde{v}(\tau,\omega))}{K_1} \tilde{y}+ \tilde{B}(\omega, \tau,  \tilde{v}(\tau,\omega)) \frac{\tilde{u}}{{K_1}})d\tau+ \frac{\tilde{Z}}{\sqrt{K_1}} d\tilde{W}_\tau
\end{equation} and we augment it with the initial and final conditions
\begin{equation}
\label{sol-S}
 \tilde{y}(0,\omega) = y^0, \ \ E(\tilde{y}(K_1  T,\cdot)) = y^T.
\end{equation} 
 Next, for $\omega \in \Omega$, we introduce an operator $\cT(\omega):\cC([0,T])\to \cC([0,T])$ which maps $\tilde{v}(\cdot,\omega)$ to a solution 
 $\tilde{y}(\cdot,\omega)$ of \eqref{stoch-K}, \eqref{sol-S} which is given by
 \begin{align} \label{Ssolll}
  \tilde{y}(\tau,\omega) = &\exp\Big(-\int_0^\tau \frac{\tilde{A}(\omega, s,  \tilde{v}(s,\omega))}{K_1} \,ds\Big) \times \\
 & \times \Big[
 \int_0^\tau \exp \Big(\int_0^s \frac{\tilde{A}(\omega, \theta,  \tilde{v}(\theta,\omega))}{K_1} \,d\theta\Big)\left( \tilde{B}(\omega,\tau,\tilde{v}(\omega,\tau)) \frac{\tilde{u}}{{K_1}}\,d\tau+\frac{\tilde{Z}}{\sqrt{K_1}}d\tilde{W}_\tau\right)
 + y^0\Big].
 \nonumber
 \end{align} We denote
 $$
 {\cal T}_\omega(\tilde{v})=\tilde{y}
 $$ and we aim to prove the existence of the fixed point for the mapping ${\cal T}_\omega$ almost surely (i.e. for ${\bf P}$-almost every $\omega\in \Omega$).
  
  In order to estimate expression \eqref{Ssolll}, we need to bound $\tilde{u}$. If we notice that \eqref{***} are Euler-Lagrange equations corresponding to the minimization problem \eqref{min}, we see that $u=B^*\varphi$, where $\varphi$ is the solution to \eqref{BSE}. This is a BSDE and we can use the standard procedure to get the necessary estimate (see \cite{Mao}). Indeed, from the It\^o formula for $|\varphi(\cdot)|^2$, we have
 $$
 \frac{1}{2}d |{\varphi}|^2=\left( \langle A^* \varphi,\varphi \rangle + \frac{|\tilde{Z}|^2}{2} \right)dt +\langle \varphi ,\tilde{Z} \rangle d{W}_t. 
 $$ By integrating the later expression on $(t,T)$ and finding the expectation of the obtained expression, while taking into  account the non-negative definitness of the matrix $A^*$, we get
 \beqs
 %\label{bound-u}
 E(|\varphi(t,\cdot)|^2)+E(\int_t^T |\tilde{Z}(t',\cdot)|^2 dt') \leq |y^T|^2.
 \eeqs Thus, since $u=B^*\varphi$ and $\varphi$ is almost surely continuous we have that the control $u$ is also almost surely continuous (and thus almost surely bounded) on $[0,T]$. Moreover, according to the results from \cite[Theorem 5.1]{7}, we know that the martingale part $Z$ of \eqref{S-Lprob} is bounded (more precisely, it is of c\'adl\'ag type). Thus, we can choose $K_1\geq 1$ such that almost surely 
 \beq \label{estSL}
 \begin{split}
 &\Big| \frac{\tilde{A}(\omega, \tau,  \tilde{v}(\tau,\omega))}{K_1} \Big| \leq 1, \ \  \frac{|B(\omega,s,\tilde{v}(\omega,s)) \tilde{u}(s,\omega)|}{K_1}\leq 1 \\ 
 &\Big|\frac{\tilde{u}}{K_1}\Big| \leq 1,  \ \ {\rm and} \ \ \Big|\frac{\tilde{Z}^2}{K_1}\Big| \leq 1.
\end{split} 
 \eeq Having this in mind, we can essentially repeat the procedure from the previous section. We square \eqref{Ssolll} and apply the expectation operator. After taking into account the It\^o isometry we get
\beqs %\label{y-0-eps-S}
 \begin{split}
E( |\tilde{y}(\tau,\omega)|^2) \leq & 3 E\Big(\exp\Big(-2\int_0^\tau \frac{\tilde{A}(\omega, s,  \tilde{v}(s,\omega))}{K_1}\,ds\Big)\times 
\\& \times\left[ \left(\int_0^\tau \exp\Big(\int_0^s \frac{\tilde{A}(\omega, \theta,  \tilde{v}(\theta,\omega))}{K_1} \,d\theta\Big) \frac{|B(\omega,s,\tilde{v}(\omega,s)) \tilde{u}(s,\omega)|}{K_1}  \,ds \right)^2\right.\\
&\qquad\left. + \int_0^\tau \exp\Big(\int_0^s 2 \frac{\tilde{A}(\omega, \tau,  \tilde{v}(\tau,\omega))}{K_1} \,d\theta\Big) \frac{|\tilde{Z}(s,\omega)|^2}{K_1}  \,ds + |y^0|^2 \Big)\right]\\& \qquad\qquad
\leq \bar{C} (\tau^2+\tau+E|y^0|^2),
\end{split}
 \eeqs for a constant $\bar{C}>0$. Next, we need to estimate the $\tau$-variation of  $E(|\tilde{y}(\tau, \omega))$.

To this end, we rewrite \eqref{stoch-K} in the integral form and again combine the expectation operator with the It\^o isometry. We have for any $0<\Delta t<1$
 \begin{equation*} 
 %\label{Sestimate}
 \begin{split}
& E\left(\Big| \tilde{y}(\tau+\Delta \tau)- \tilde{y}(\tau) \Big|^2\right)\\& \leq  2 E \left[\left(\int_\tau^{\tau+\Delta \tau} \left(-\frac{\tilde{A}(\omega, s,  \tilde{v}(s,\omega))}{K_1}\tilde{y}(s,\cdot)+ \frac{B(\omega,s,\tilde{v}(\omega,s))\tilde{u}(s,\cdot)}{K_1}\right) ds \right)^2 \right]\\&+2 E \left[ \int_\tau^{\tau+\Delta \tau} \frac{\tilde{Z}(s,\cdot)^2}{K_1} ds \right].
 \end{split}
 \end{equation*} From here and \eqref{estSL}, we have
$$
 E\left(\frac{\Big| \tilde{y}(\tau+\Delta \tau)- \tilde{y}(\tau) \Big|^2}{\Delta \tau}\right) \leq C.
$$  Hence, for any $\tilde{T}>0$,  $|\tilde{y}|^2$ and its H\"older norm $ \|\tilde{y} \|_{C^{1/2}([0,T])}$ are uniformly bounded for ${\bf P}$-almost every $\omega\in \Omega$ on $[0,\tilde{T}]$. 
 
 Thus, we can repeat the procedure from the proof of Theorem \ref{main-T} (after relation \eqref{estimate}) to conclude  the existence of the fixed point for the mapping $\cT_\omega$ for ${\bf P}$-almost every $\omega\in \Omega$. The obtained function $y$ is the solution to \eqref{S-Lprob}, \eqref{L-in-fin}.
\ep

\section{An application on population dynamics}

A motivation for considering system \eqref{porous-med} was the porous media equation which raises a constant interest from the viewpoint of both purely theoretical or applied mathematics. It has the form
\begin{equation}
\label{PM}
\pa_t u=\nabla \big( (m-1) |u|^{m-1} \nabla u \big)
\end{equation} and it describes the standard model of gas flow
through a porous medium (Darcy-Leibenzon-Muskat), nonlinear
heat transfer (Zel'dovich-Raizer), Boussinesq's groundwater flow, population dynamics (Gurtin-McCamy) etc. \cite{V}. Different variant of the equation such as stochastic \cite{BPR} or fractional \cite{PQRV} are also well known and well investigated.

\begin{figure}[h!]
  \centering
  \begin{subfigure}[b]{0.4\linewidth}
    \includegraphics[width=\linewidth]{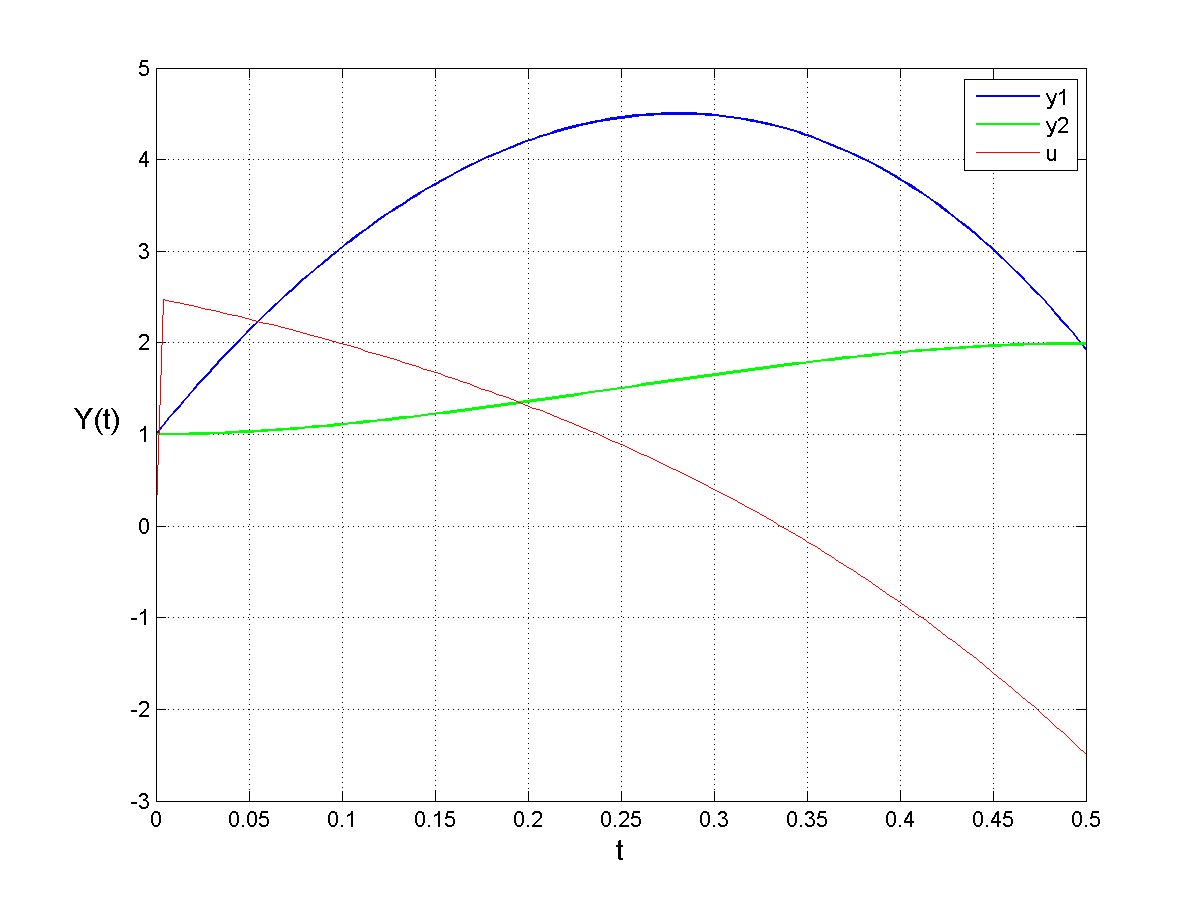}
    \caption{Initial iteration}
  \end{subfigure}
  \begin{subfigure}[b]{0.4\linewidth}
    \includegraphics[width=\linewidth]{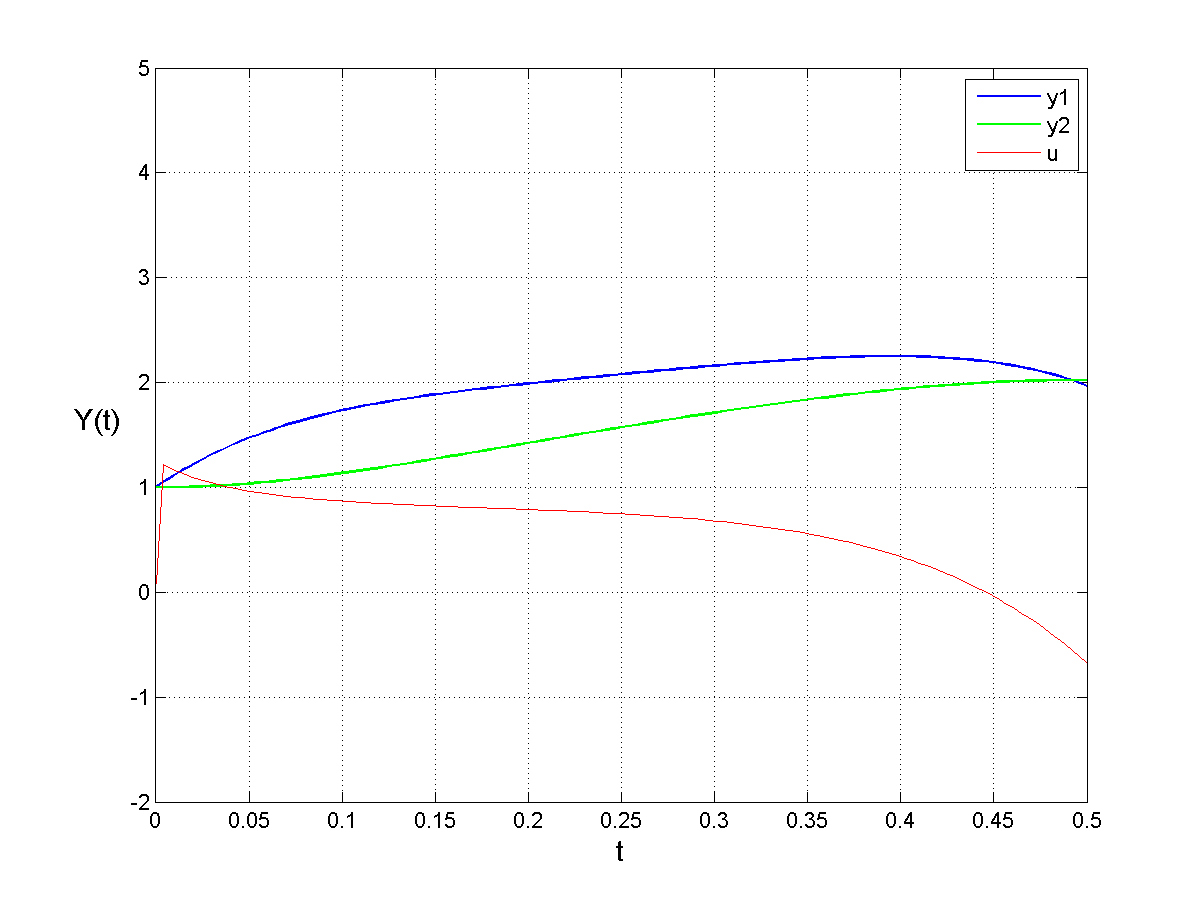}
    \caption{After two iterations.}
  \end{subfigure}
  \begin{subfigure}[b]{0.4\linewidth}
    \includegraphics[width=\linewidth]{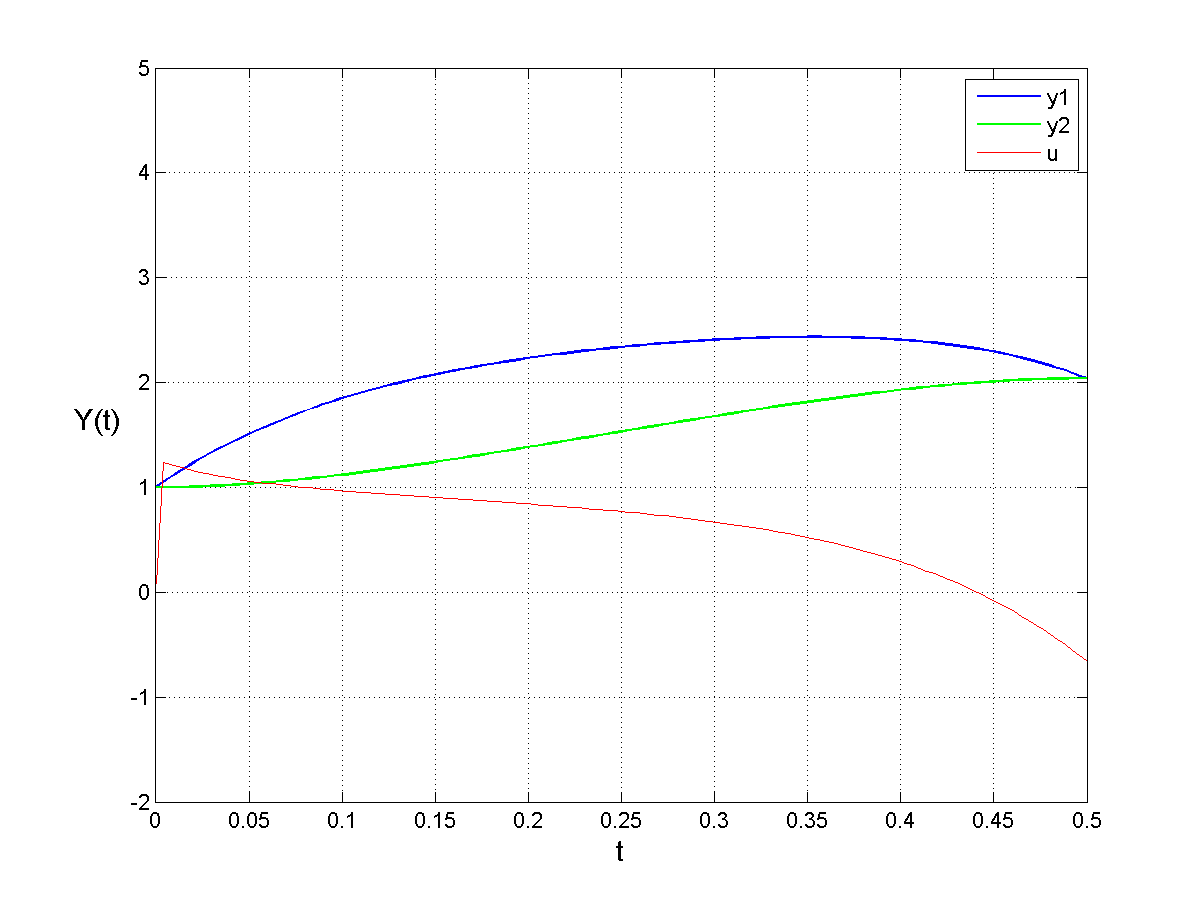}
    \caption{After three iterations.}
  \end{subfigure}
  \begin{subfigure}[b]{0.4\linewidth}
    \includegraphics[width=\linewidth]{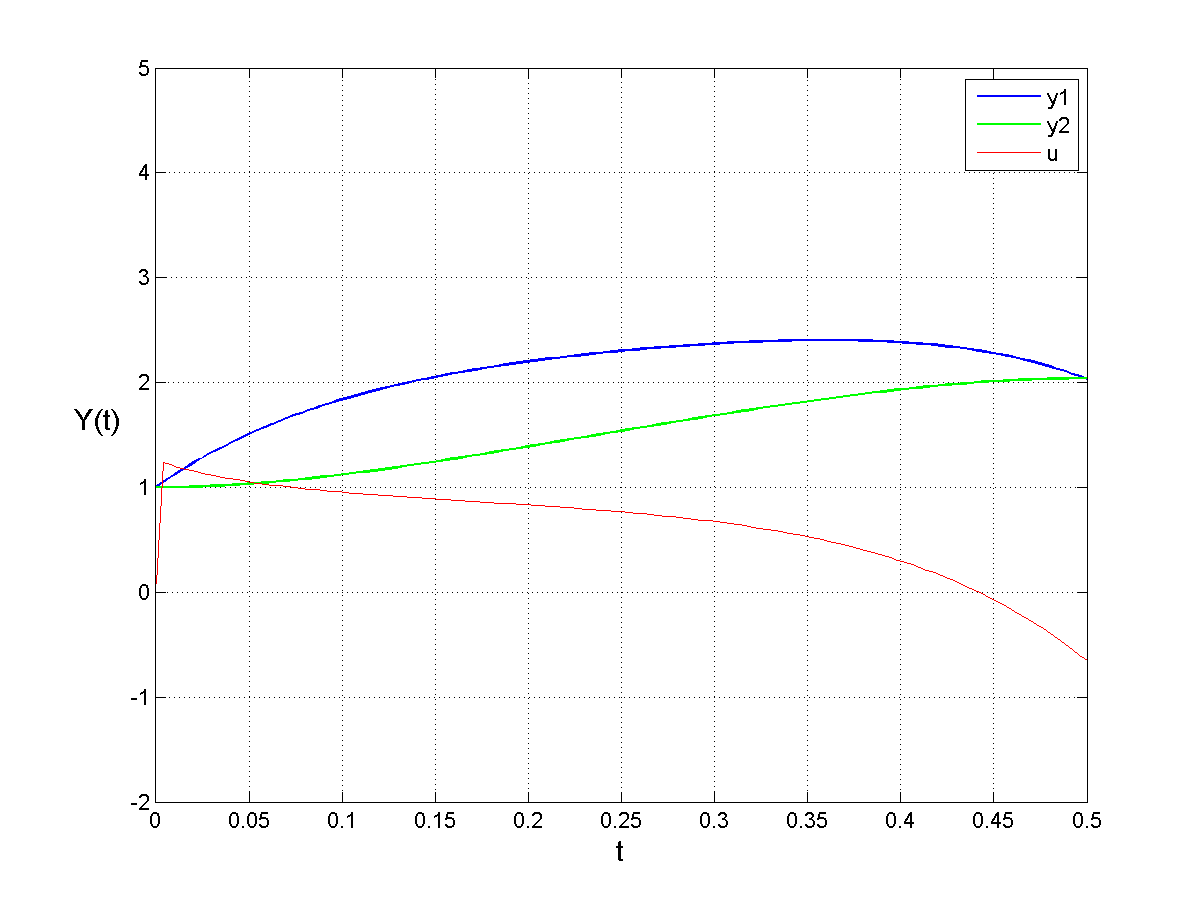}
    \caption{After four iterations.}
  \end{subfigure}
  \begin{subfigure}[b]{0.4\linewidth}
    \includegraphics[width=\linewidth]{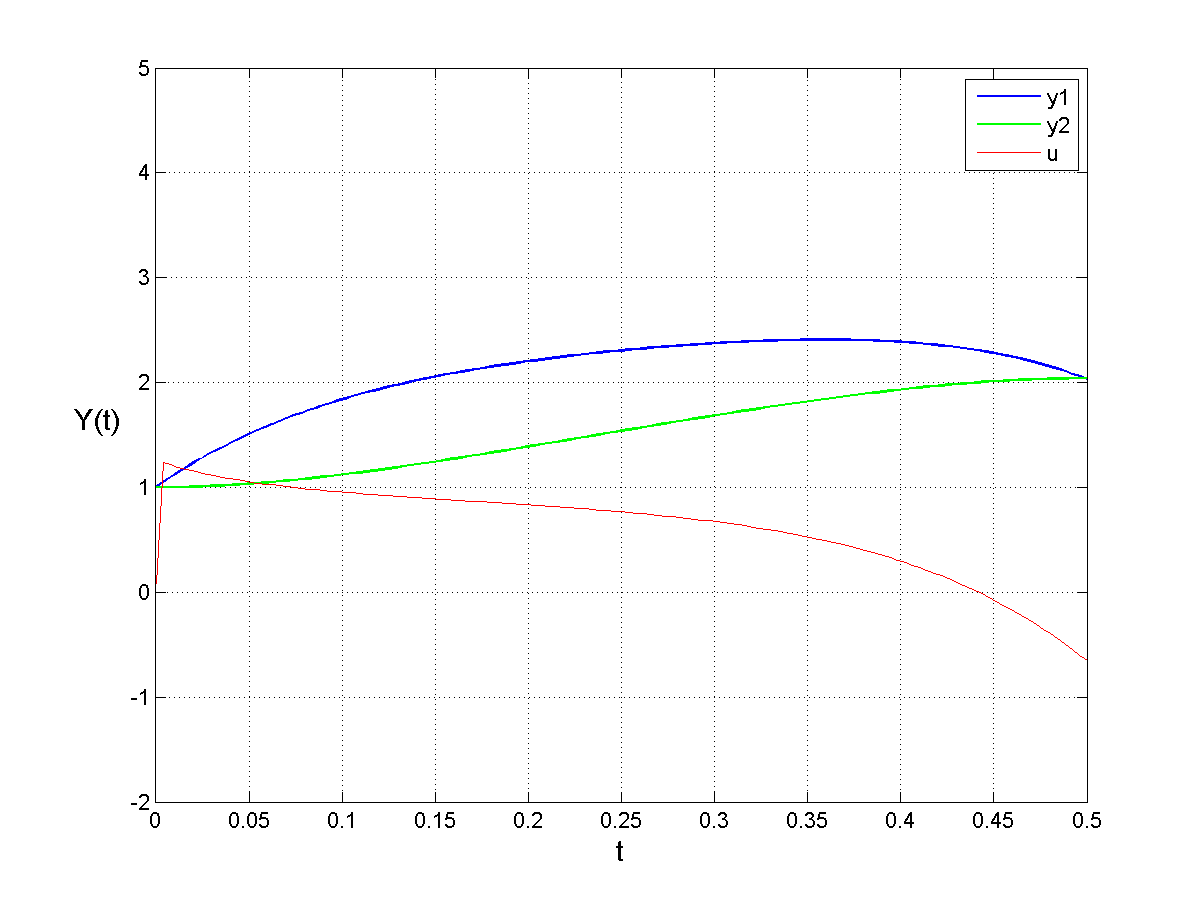}
    \caption{After five iterations.}
  \end{subfigure}
  \begin{subfigure}[b]{0.4\linewidth}
    \includegraphics[width=\linewidth]{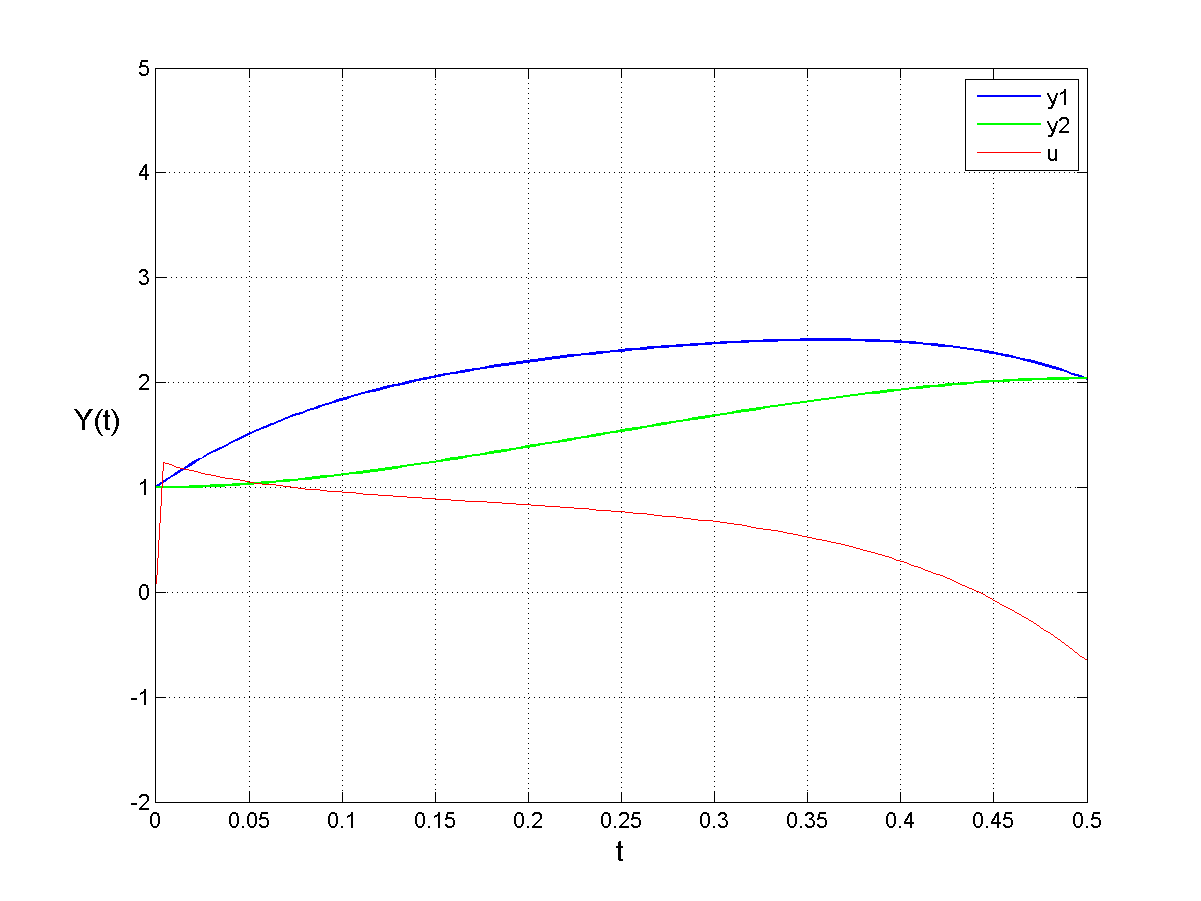}
    \caption{After six iterations.}
  \end{subfigure}  
  \caption{The numerical procedure after one to six iterations involving perturbations of final states.}
  \label{fig:2}
\end{figure}

\begin{figure}[h!]
  \centering
  \begin{subfigure}[b]{0.4\linewidth}
    \includegraphics[width=\linewidth]{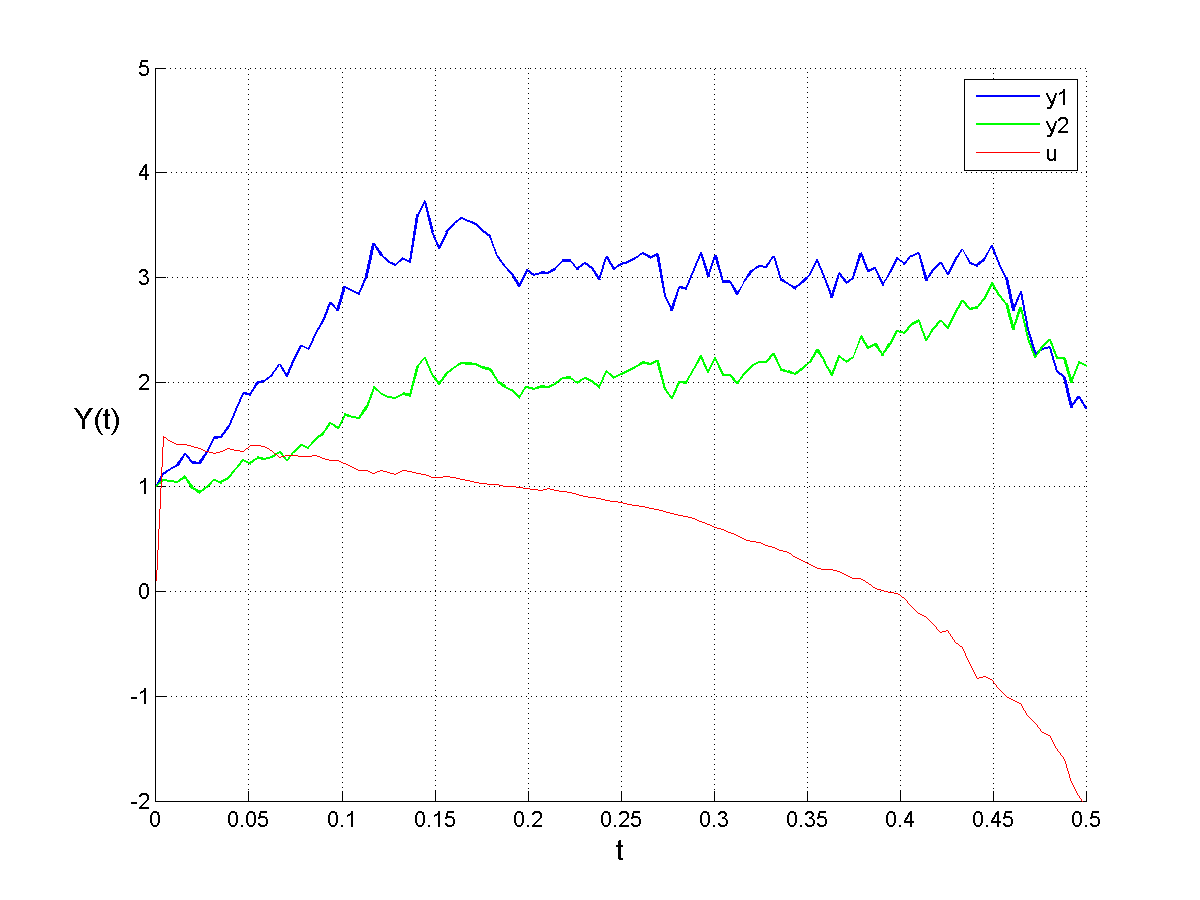}
    %\caption{}
  \end{subfigure}
  \begin{subfigure}[b]{0.4\linewidth}
    \includegraphics[width=\linewidth]{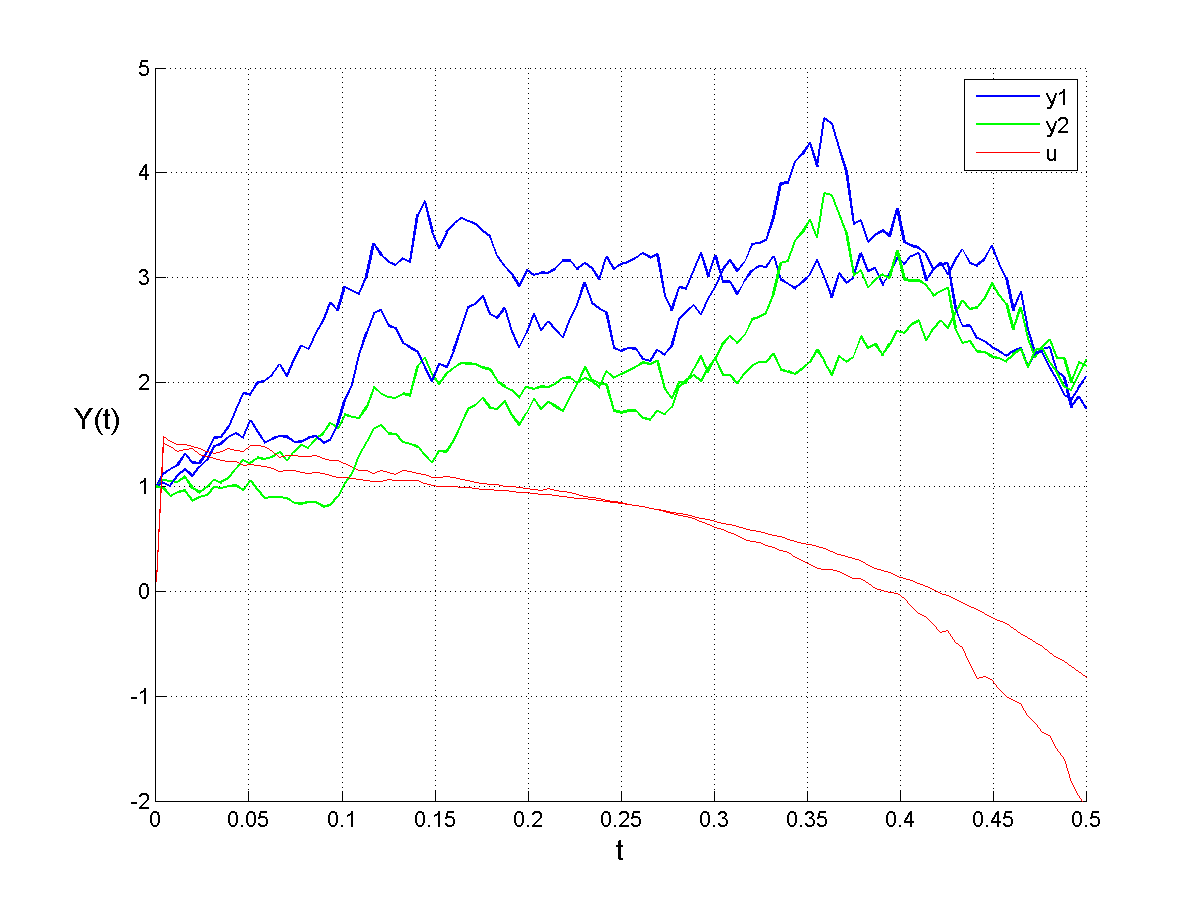}
    %\caption{More coffee.}
  \end{subfigure}
  \begin{subfigure}[b]{0.4\linewidth}
    \includegraphics[width=\linewidth]{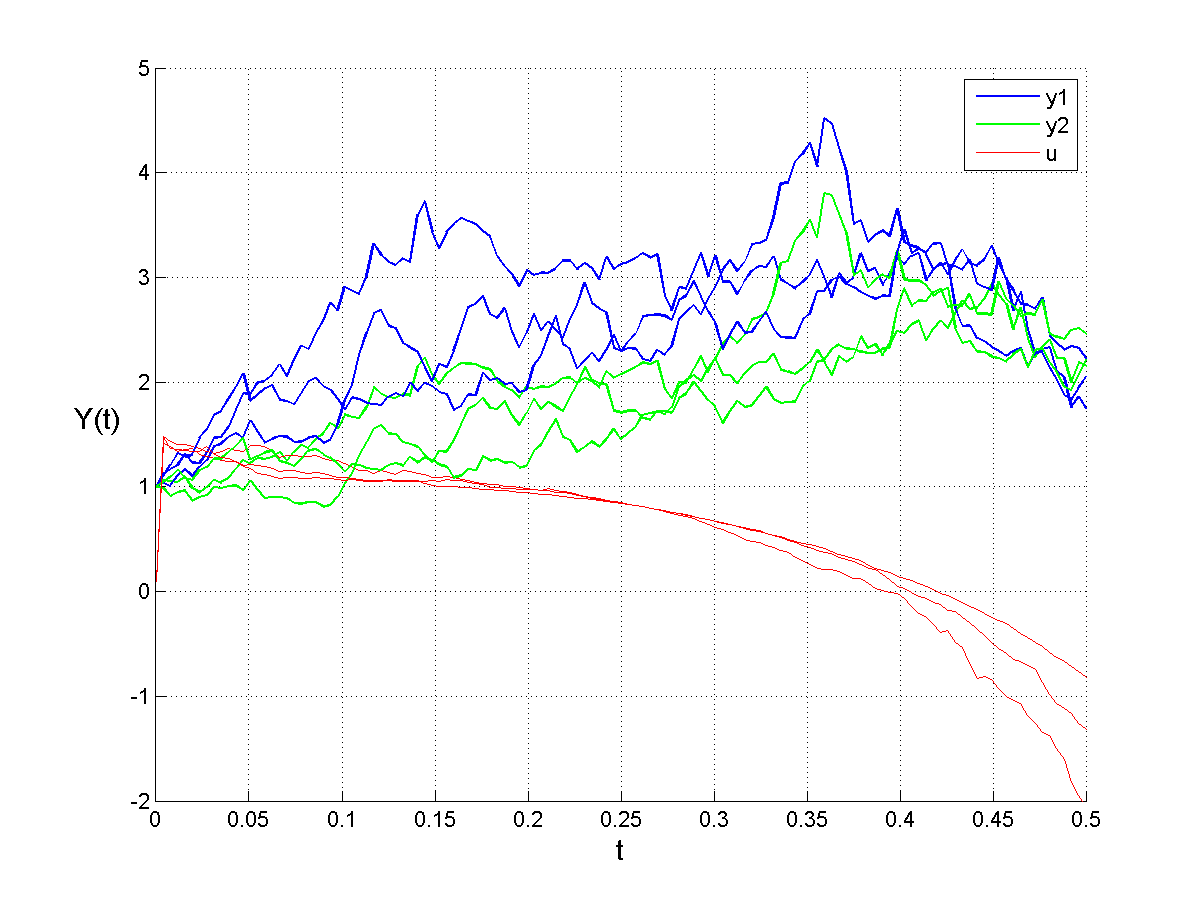}
    %\caption{More coffee.}
  \end{subfigure}
  \begin{subfigure}[b]{0.4\linewidth}
    \includegraphics[width=\linewidth]{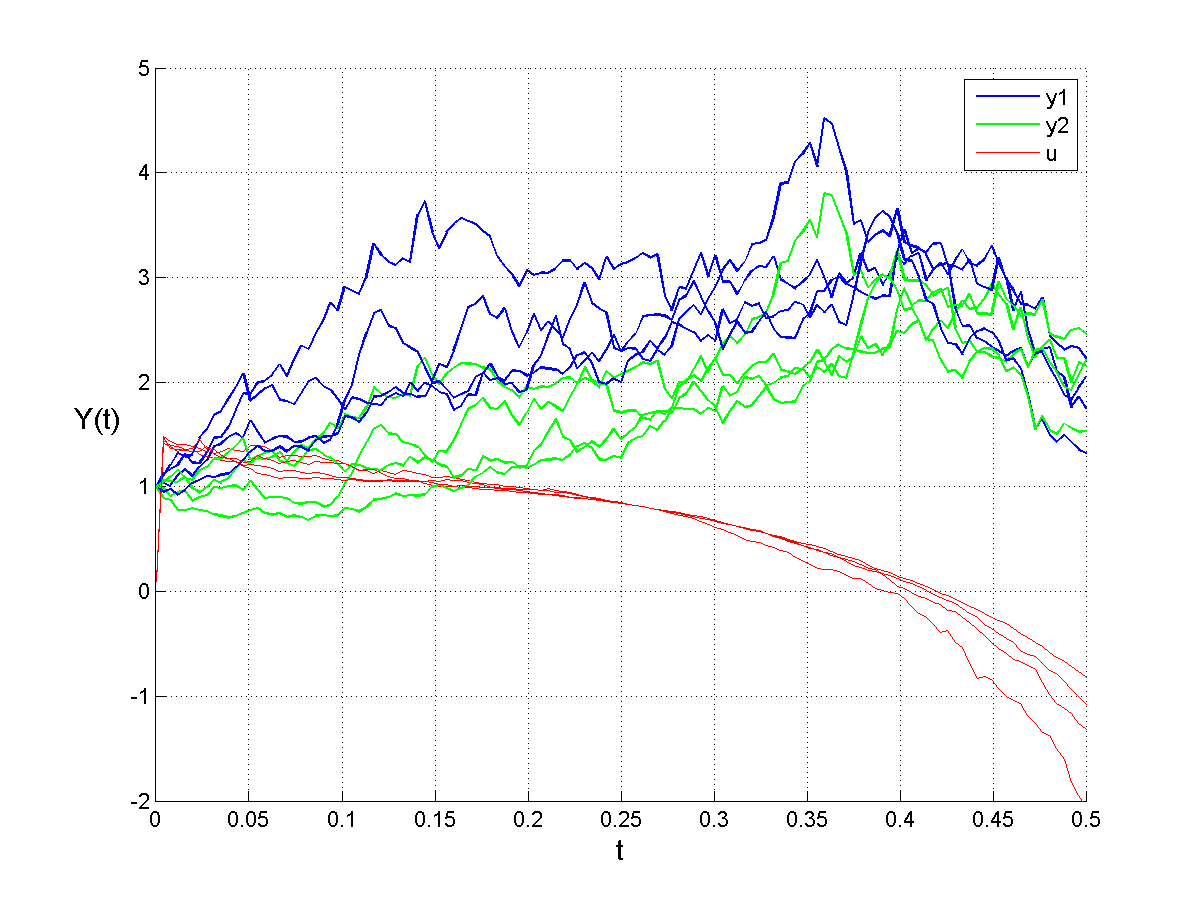}
    %\caption{More coffee.}
  \end{subfigure}
  \begin{subfigure}[b]{0.4\linewidth}
    \includegraphics[width=\linewidth]{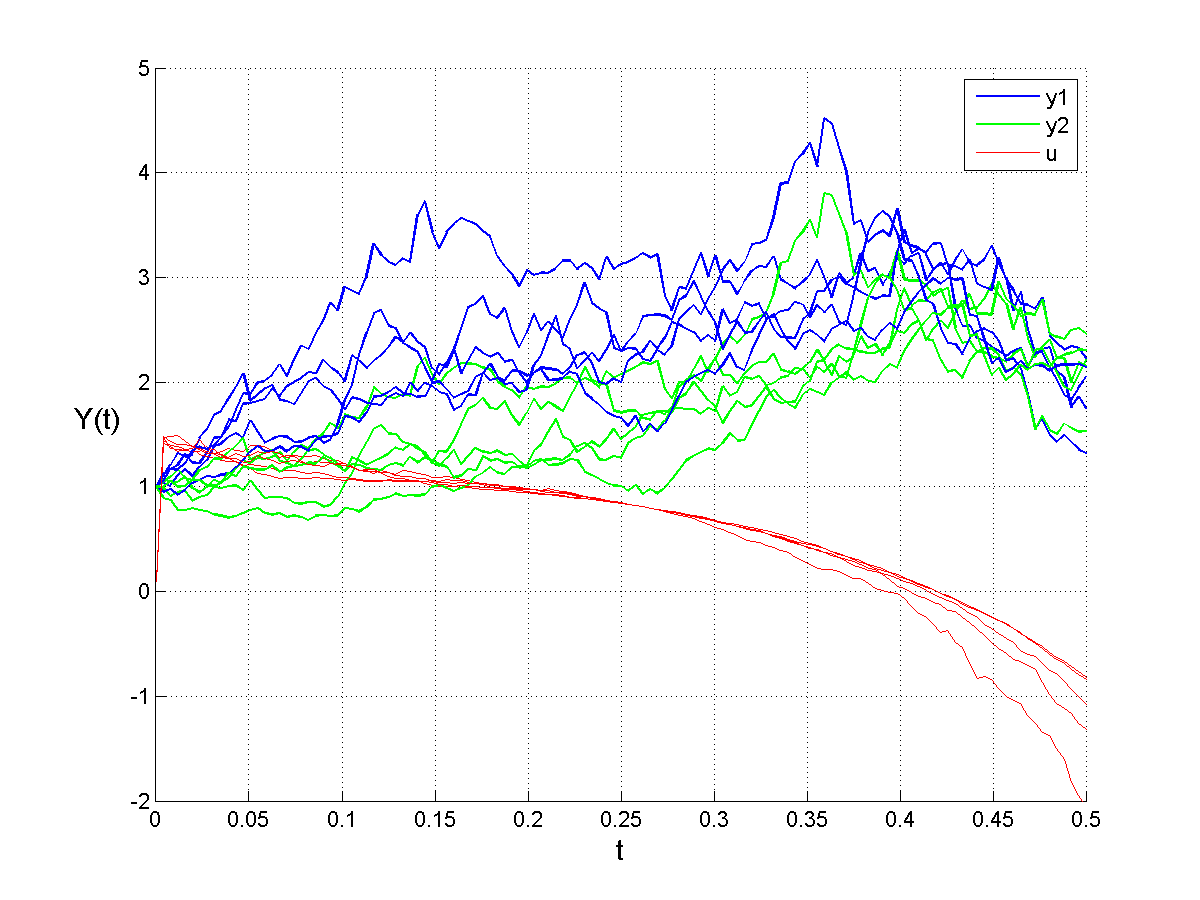}
    %\caption{More coffee.}
  \end{subfigure}
  \begin{subfigure}[b]{0.4\linewidth}
    \includegraphics[width=\linewidth]{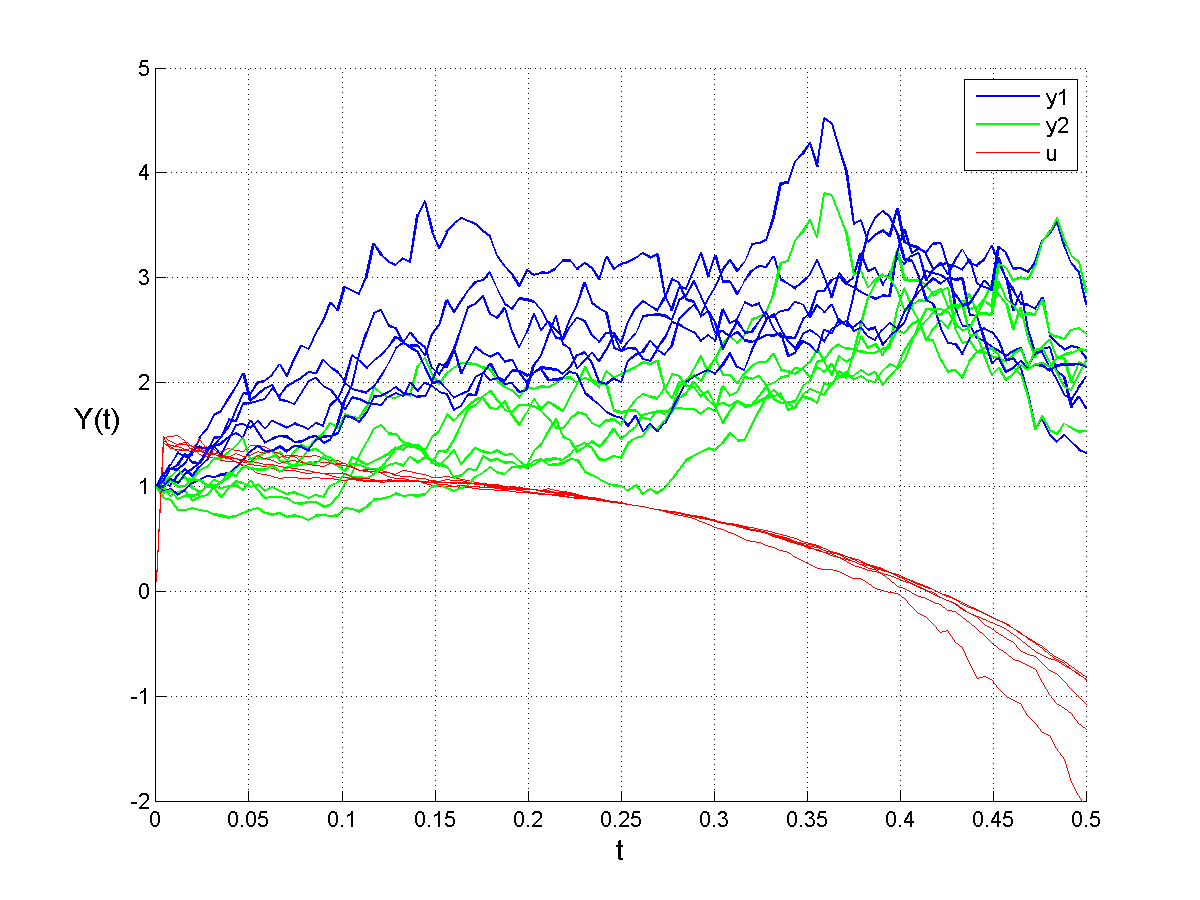}
    %\caption{More coffee.}
  \end{subfigure}  
  \caption{The state evolution for six different samplings.}
  \label{fig:3}
\end{figure}

To this end, we shall consider an ODE counter-part of \eqref{PM}. More precisely, we shall consider population dynamics with species intending to avoid crowding. It is interesting to inspect how to control the populations by adding new individuals or in some other way by improving the living conditions (in the frame of the given nonlinear model, of course). We assume that we have two species whose population quantities are denoted by $y_1$ and $y_2$ and which have tendency to avoid crowding. Mathematical model of the phenomenon is given by
\beq
\label{example}
\begin{split}
dy_1&=(|y_1+y_2|(-2 y_1+2y_2) + u) dt+Z_1 dW_t\\
dy_2&=|y_1+y_2| (y_1-y_2) dt+Z_2 dW_t
\end{split}
\eeq and we aim to maintain the population by randomly introducing new individuals $u_1$ and $u_2$ of the corresponding species into the system. Thus, we start with the populations
$$
y_1(0,\omega)=y_2(0,\omega)=1,
$$ and we want to have the same population at a final time $T=0.5$
\begin{equation}
\label{fin-num}
E(y_1(0.5,\cdot))=E(y_2(0.5,\cdot))=2.
\end{equation} By direct substitution, we know that for the problem
$$
dy=(A(t)y+ B(t)u)dt+ZdW_t, \ \ E(y(0,\cdot))=y_0, \ \ E(y(0.5,\cdot))=y_1, 
$$ the control function is given by 
\begin{equation}
\label{control-expl}
u(t,\omega)=-B^{\rm T} \Phi^{\rm T}(0,t) W_c(0,1)^{-1}[y_0-\Phi(0,1)y_1]
\end{equation} where 
\begin{equation}
\label{wc-int}
W_c(0,t)=\int_0^t \Phi(0,\tau)B(\tau)B^{\rm T}(\tau)\Phi^{\rm T}(0,\tau) d\tau
\end{equation} and $\Phi$ solves the system
\begin{equation}
\label{id-id}
\Phi'=A(t) \Phi, \ \ \Phi(0,0)=I.
\end{equation}
In order to solve \eqref{example}, we use the following recursive procedure
\beqs
%\label{example-1}
\begin{split}
dy^n_1&=(|y^{n-1}_1+y^{n-1}_2|(-2 y^n_1+2y^n_2) + u^n) dt+Z^n_1 dW_t\\
dy^n_2&=|y^{n-1}_1+y^{n-1}_2| (y^n_1-y^n_2 )dt+Z^n_2 dW_t
\end{split} 
\eeqs where $u^n$ is given by the corresponding variant of \eqref{control-expl}. The results are shown on Figures \ref{fig:2} and  \ref{fig:3}. Please note that instead of the control $u$, we plot the scaled control $u/10$ in order to provide better presentation of the results.

\begin{figure}[h!]
  \centering
  \begin{subfigure}[b]{0.4\linewidth}
    \includegraphics[width=\linewidth]{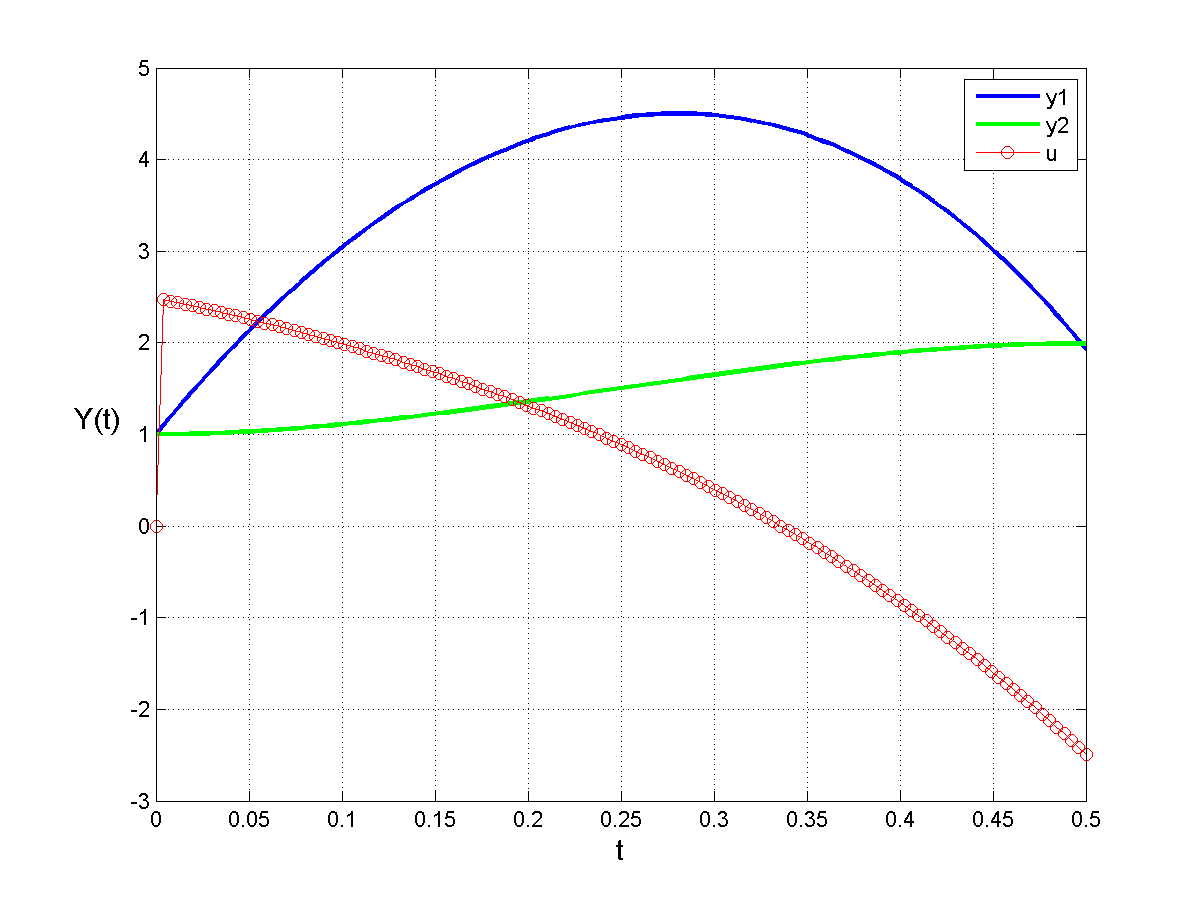}
    \caption{Initial iteration}
  \end{subfigure}
  \begin{subfigure}[b]{0.4\linewidth}
    \includegraphics[width=\linewidth]{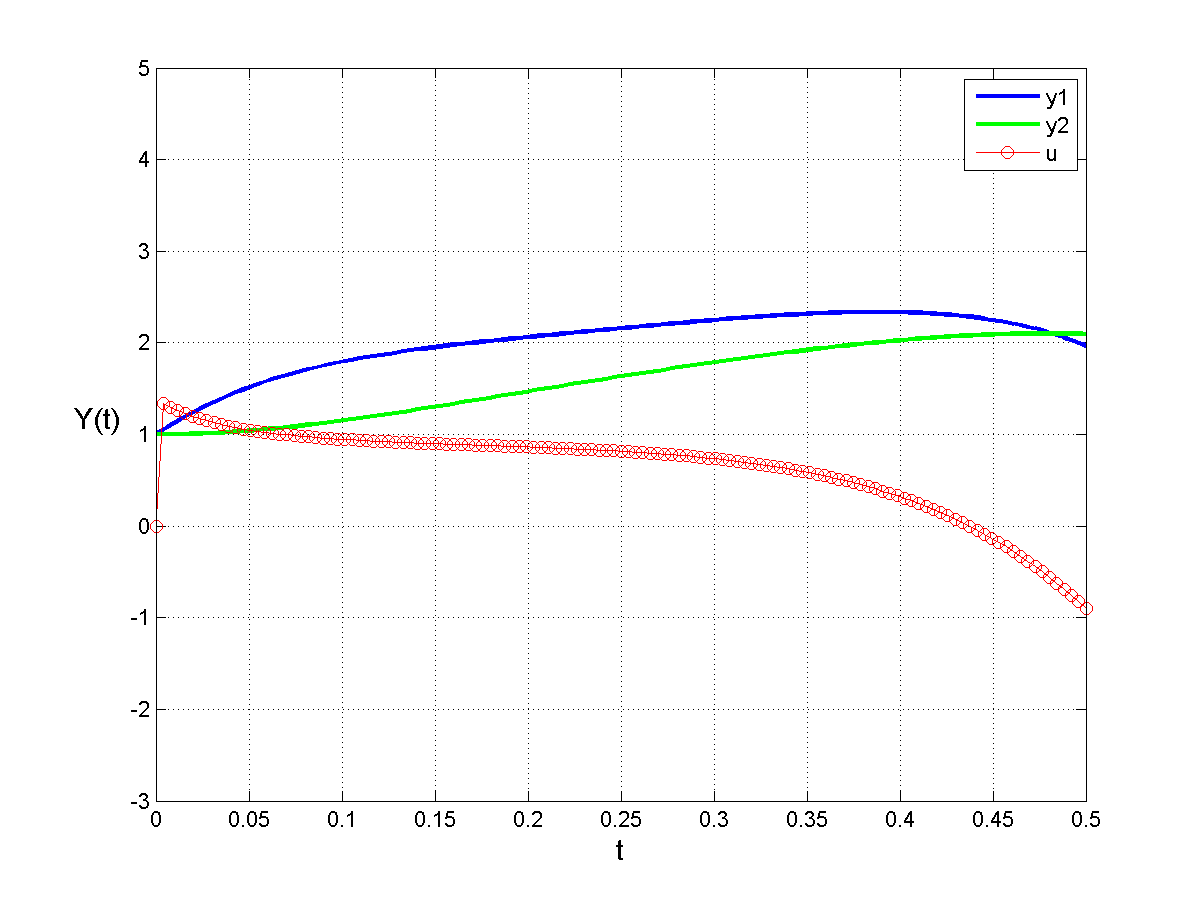}
    \caption{After two iterations.}
  \end{subfigure}
  \begin{subfigure}[b]{0.4\linewidth}
    \includegraphics[width=\linewidth]{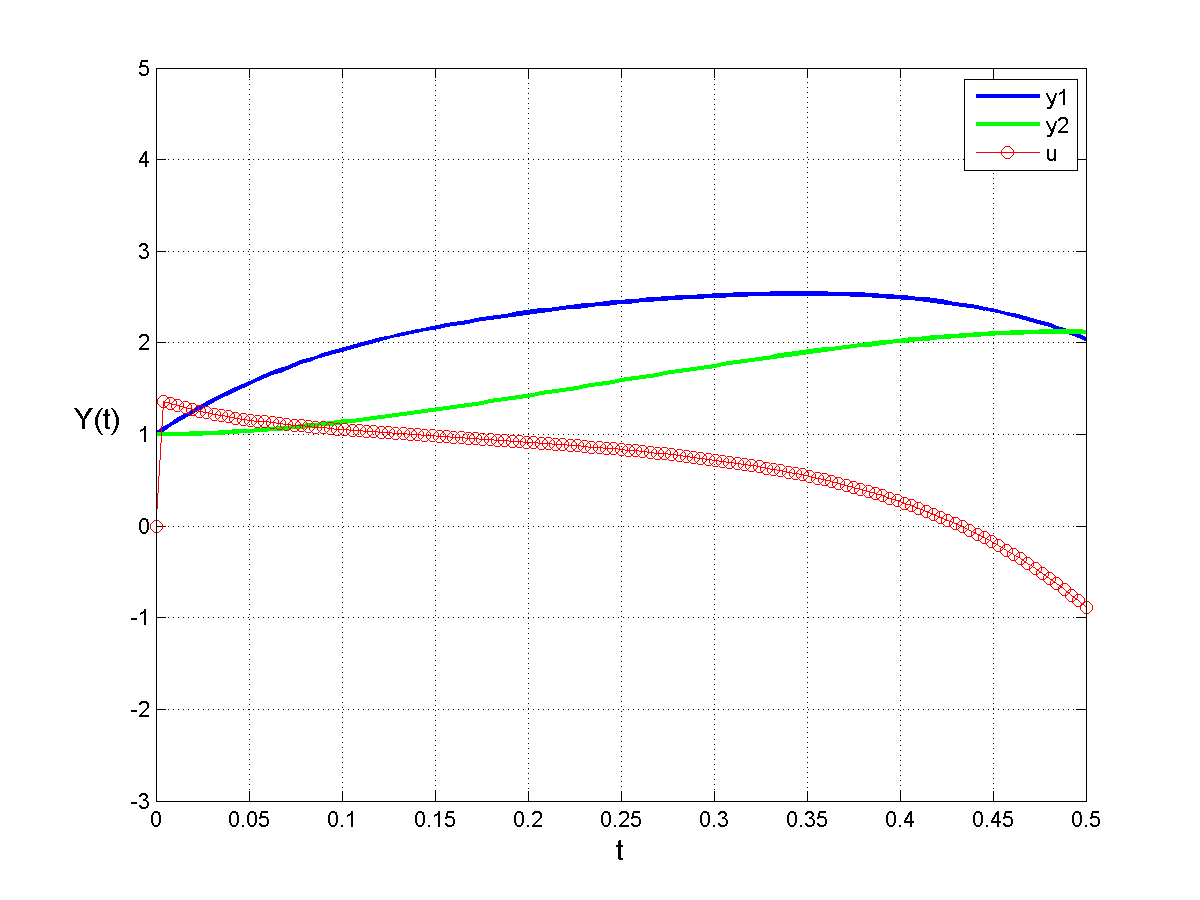}
    \caption{After three iterations.}
  \end{subfigure}
  \begin{subfigure}[b]{0.4\linewidth}
    \includegraphics[width=\linewidth]{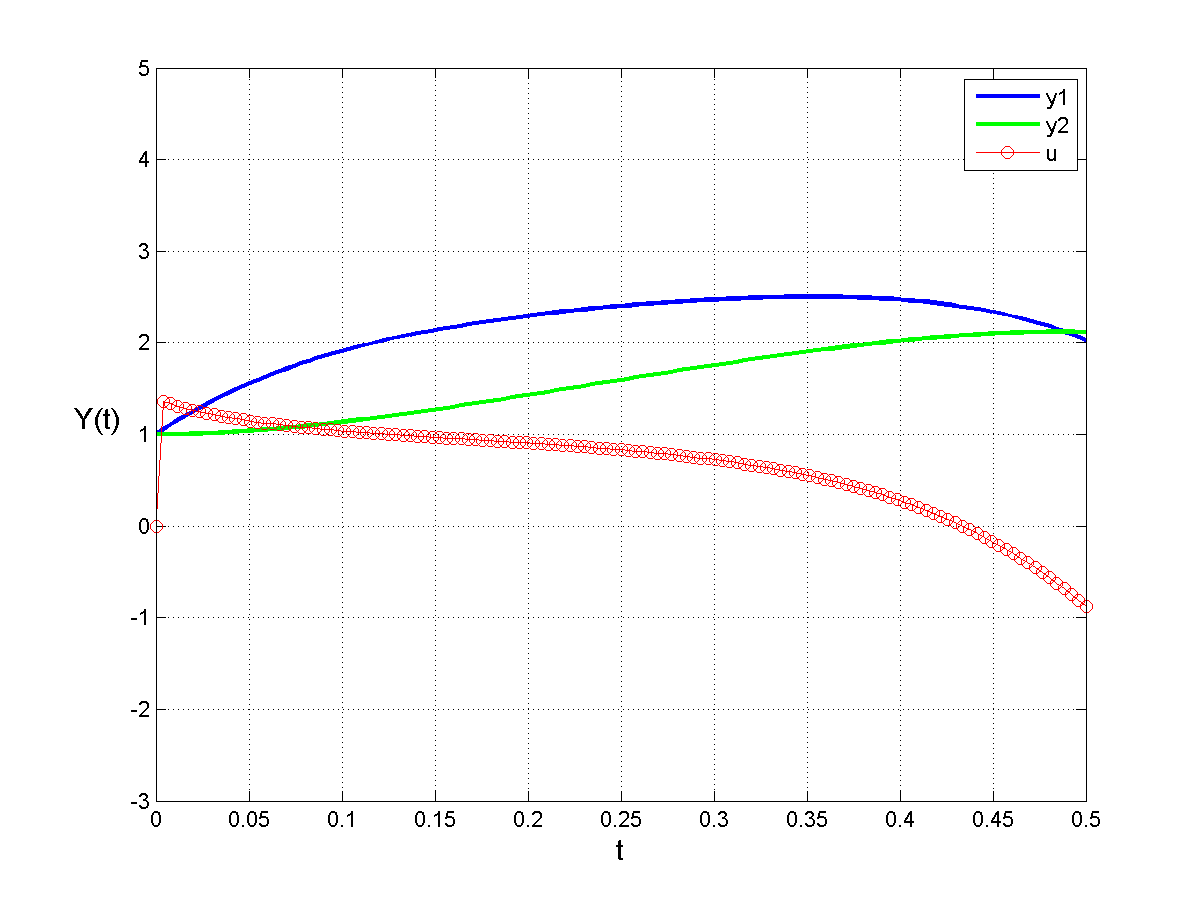}
    \caption{After four iterations.}
  \end{subfigure}
  \begin{subfigure}[b]{0.4\linewidth}
    \includegraphics[width=\linewidth]{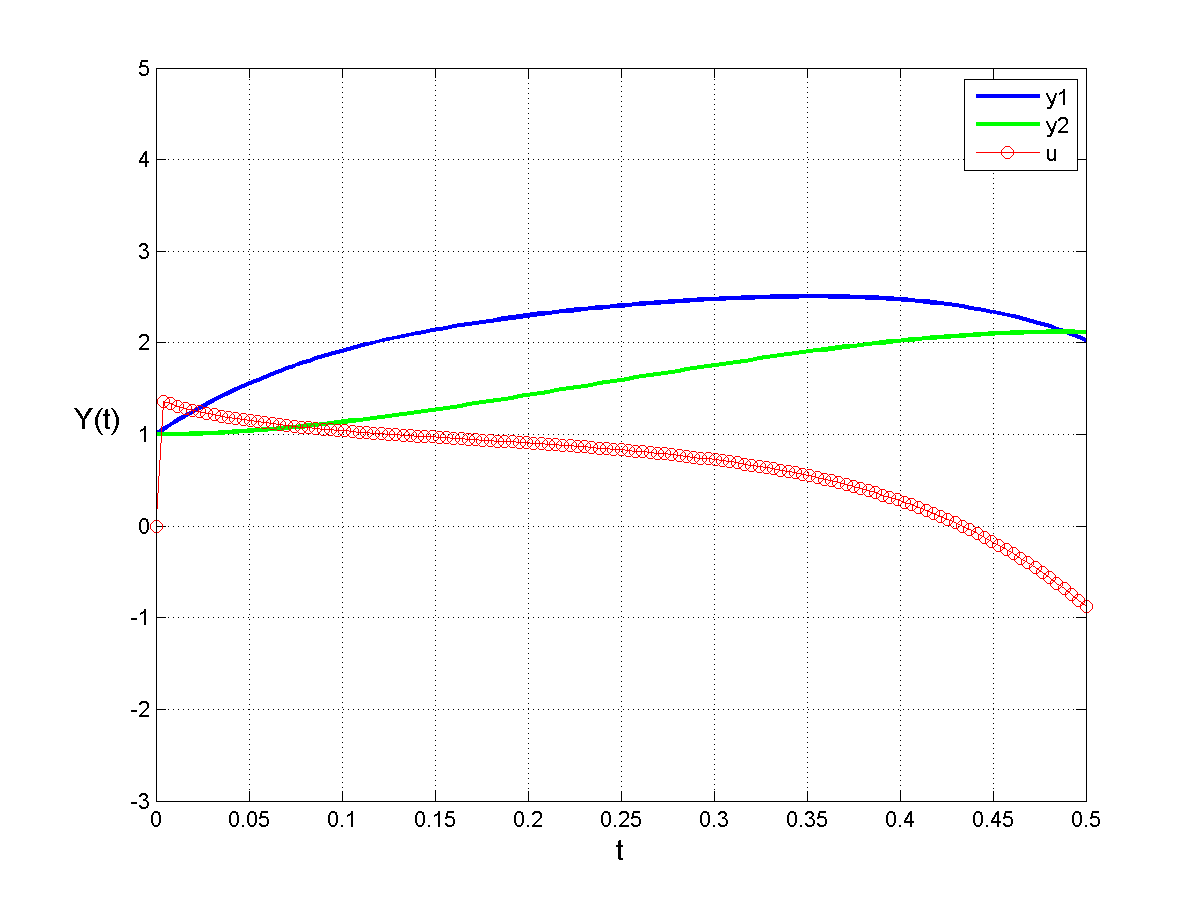}
    \caption{After five iterations.}
  \end{subfigure}
  \begin{subfigure}[b]{0.4\linewidth}
    \includegraphics[width=\linewidth]{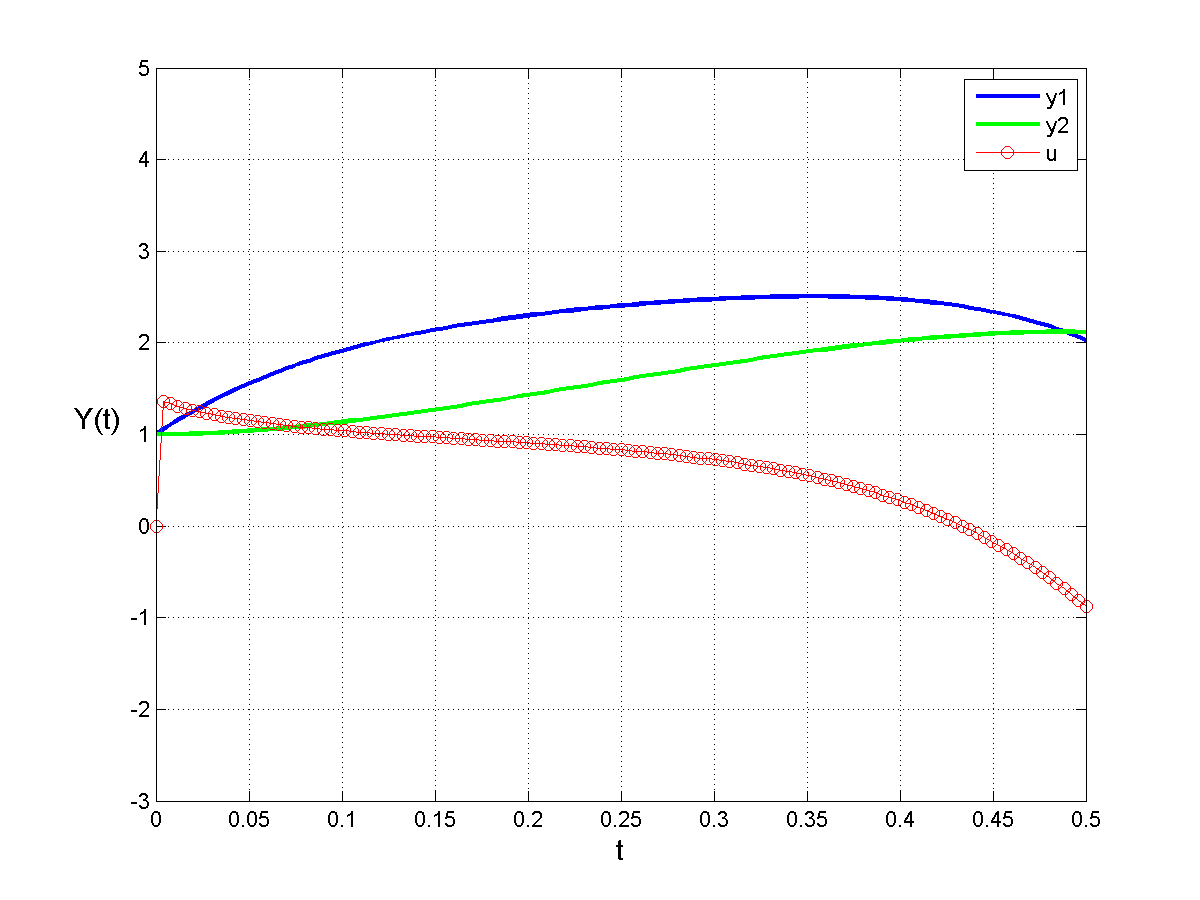}
    \caption{After six iterations.}
  \end{subfigure}  
  \caption{The numerical procedure after one to six iterations.}
  \label{fig:1}
\end{figure}

As we can see on Figure \ref{fig:1} showing the control $u$ and the approximate state $(y_1,y_2)$ after $n=1,2,3,4,5, 6$ iterations in the deterministic situation, our method provides the exact control for the given nonlinear system. We notice that the procedure stabilizes already after the third iteration, for instance, the relative change of $u$ from fourth to fifth iteration is less than $1\%$ in the $L^\infty$-norm.

%%%%%%%%%%%%%%

We can also notice that the final states are not exactly at the prescribed values $y_1(0.5)=y_2(0.5)=2$ which is due to the numerous computations of various integrals \eqref{wc-int}  and matrix exponential \eqref{id-id}. More precisely, we need to compute $\mathcal{O}(1/dt)$ integrals in every iteration until the fixed point is reached, so the appearance of the numerical error is not unexpected.  For instance, approximation of the matrix exponential by the polynomial of the second order gives even worse results since in this case we end up in 
$$
y_1(0.5) = 1.7640, \ \ y_2(0.5) = 2.1919.
$$
On the other hand, approximations of $e^A$ by the polynomials of fifth and sixth order provide the final states
\begin{align*}
y_1(0.5) = 1.9851, \ \ y_2(0.5) = 2.1195 \ \ \text{ fifth order approximation of the exponential.}\\
y_1(0.5) = 2.0209, \ \ y_2(0.5) = 2.1070 \ \ \text{ sixth order approximation of the exponential.}
\end{align*} Interestingly, the approximation of seventh order only slightly changes the final states. 

In order to improve the procedure, we introduce the perturbation of the final states after each iteration as follows: after the first iteration, the reached values of the final states are $y_1(0.5)= 1.93$ and $y_2(0.5)=1.99$ while the prescribed final state was $y_1(0.5)= y_2(0.5)=2$. For the next iteration, we shall therefore require $y_1(0.5)= 2 +  0.07\alpha$ and $y_2(0.5)=2 + 0.01\alpha$,  where $\alpha$ is the relaxation parameter and in our experiments $\alpha = 10^{-1}$. In this way, we expect to compensate for the numerical errors by purposely introducing a perturbation of the final states. 
Given that we have omitted the stochastic forcing, we are able to fine tune the perturbation and obtain almost perfect results as demonstrated on Figure \ref{fig:2}.

%%%%%%%%%%%%%%

An extension of this approach to the stochastic case is shown in Figure \ref{fig:3}, where we have a simulation after five iterations for six different samplings. We see that the final state is not reached exactly as the stochastic term distorts evolution of the state functions from the initial to the final conditions. Instead, the final state is reached in an averaged sense \eqref{fin-num}. In order to demonstrate this, we have performed five sets of numerical experiment with twenty different samplings of the stochastic forcing in each set, while limiting the maximal number of iterations to $n = 1$ in the first set of experiments, $n = 2$ in the second set, and continuing in a similar fashion to $n = 5$ for the last set of twenty experiments.  
The expected values of the final states for each set of twenty experiments as well as the total expected value (for all 100 experiments) are provided bellow in the Table \ref{tab:table1}. The total expected value is the average of final states obtained for all iterations. Although the total expected value is related to a slightly different processes, they solve similar equations and the convergence of the method seems to be fast.

\begin{table}[h!]
  \begin{center}
    \caption{Expected values after $n =1, 2,\ldots, 5$ iterations (with twenty experiments for each value of $n$) and total expected value.  {\em Exp} and {\em Iter} denote expected value and number of iterations, respectively.}
    \label{tab:table1}
    \begin{tabular}{l|c|c|c|c|c|c|} % <-- Alignments: 1st column left, 2nd middle and 3rd right, with vertical lines in between
      \textbf{Exp$\setminus$ Iter } & \textbf{1 iter} & \textbf{2 iter} & \textbf{3 iter} & \textbf{4 iter} & \textbf{5 iter} & \textbf{Total}\\
            \hline     
      $E(y_1(1,\cdot))$ & 1.8126 & 2.0025 &  1.8742 & 2.0805 & 2.1851 & {\bf 1.9910} \\
      $E(y_2(1,\cdot))$ & 1.8391 & 2.0206 &  1.9112 & 2.0942 & 2.1912 & {\bf 2.0113} 
    \end{tabular}
  \end{center}
\end{table}

\begin{rem} 
The introduced method opens new possibilities in inspecting controllability properties of the nonlinear systems. An obvious generalization of the given procedure can be performed on the systems of the form
\begin{eqnarray}
 && dy(t)= \left( - A(t,\omega,y(t,\omega),Z(t,\omega),u(t,\omega)) y(t,\omega)+B(t,\omega) u(t,\omega)+ C(t,\omega) Z(t,\omega)\right] dt \nonumber\\
 && \phantom{ dy(t)=}+\left[ D(t,\omega) y(t,\omega)+E(t,\omega) u(t,\omega)+F(t,\omega) Z(t,\omega) \right]dW(t), \nonumber \\
 &&\phantom{ dy(t)} y(0,\omega)=y^0,\qquad Ey(T,\omega)=y^T, \label{genlast}
\end{eqnarray} where $A$ is non-negative definite matrix while $B,C,D,E,F$ are well defined matrices. Combining the approach from \cite{4C} with the method introduced here, all statements from Section 3 should apply to this case which opens the possibility to solve the averaged control problem for \eqref{genlast}.
  \end{rem}
%  Regarding the stochastic case, all statements and proofs are given for the case of $L^2$ norm. It should be noted the results can be extended on th case $L^p$, for $p>2$.

%%%%%%%
 \section*{Acknowledgement}
%%%%%%%
This work is partially supported by Ministry of Education, Science and Technological Development of the Republic of Serbia (Grant No. 451-03-9/2021-14/200125). It is also supported by the Lise Meitner project number M 2669-N32 of the Austrian Science Fund (FWF) and by the Croatian Science Foundation under Project MiTPDE (number IP-2018-01-2449).
This article is based upon work from COST Action CA15225 FRACTIONAL and CA15125 DENORMS supported by COST (European Cooperation in Science and Technology). 

The permanent address of D.M. is University of Montenegro, Montenegro.

%%%%%%%
%%%%%%%

 \bibliographystyle{plain}

\renewcommand{\baselinestretch}{1}

%%%%%%%%%%%%%%%%%
%%%%%%%%%%%%%%%%%
 \end{document}